# Climate-resilient electric vehicle charging infrastructure for sustainable cities: An interpretable causal-ensemble framework for preventive maintenance and low-carbon mobility

Cande Lian [a,b,1] , Wentao Zeng[a,c,1], Jiabin Wu [a,*], Yiming Bie [d], Wei Zhou [e]

[a] *School of Management, Foshan University, Foshan 528000, China*
[b] *School of Mathematics, Foshan University, Foshan 528000, China*
[c] *School of Mechanical and Electrical Engineering and Automation, Foshan University, Foshan 528225, China*
[d] *School of Transportation, Jilin University, Changchun 130022, China*
[e] *Department of Civil and Environmental Engineering, National University of Singapore, Singapore 117576, Singapore*
* **Corresponding author: Jiabin Wu** (E-mail address: jiabinwu@fosu.edu.cn)
[1] The two authors contribute equally

**Abstract:** Reliable electric vehicle (EV) charging infrastructure is a cornerstone of sustainable, low-carbon cities, yet urban climate stress such as extreme heat, heavy precipitation, and humidity increasingly raises equipment fault risk and undermines the resilience of urban energy and mobility services. Shifting operation from reactive repair to preventive maintenance depends on accurate, forward-looking fault-risk prediction, a task complicated by the heterogeneous time scales of physical, behavioral, contextual, and historical signals and by forecasting over a multi-week horizon. We develop FGDSE, a feature-governed dynamic stacking ensemble that forms an interpretable decision-support system for climate-resilient charging-asset management. It partitions heterogeneous signals into four feature families, assigns each to a domain expert whose inductive bias matches the data, and adds two deep temporal experts for short-term pulses and long-term degradation; a horizon-wise gating mechanism then learns adaptive weights to forecast daily fault risk over 1 to 30 days. SHAP attribution and an X-learner extend the probabilistic output into causal decision support with post-level treatment effects. On 25 months of data from 13 stations, FGDSE surpasses twelve baselines beyond the ten-day horizon, sustains about 85% macro-recall at 30 days with an AUC decay of only 3.2 points, and reveals a shift of dominance from fault history toward climate stress. It identifies extreme heat as the sole exposure whose causal effect amplifies over time, flagging roughly 30% of posts as heat-sensitive and yielding quantitative thresholds for climate-adaptive maintenance that strengthens urban mobility resilience and sustains low-carbon travel.

# 1. Introduction

Cities are the central arena of the global transition to low-carbon, climate-resilient development, and decarbonizing urban mobility is a decisive part of this effort. As transport electrification accelerates, public charging infrastructure is evolving from a supporting amenity into critical urban energy-mobility infrastructure that links the power grid with the transport system and underpins reliable low-carbon travel, while conventional fueling stations face growing pressure to transition toward charging services and adaptive reuse (Zhou et al., 2026). The International Energy Agency's Global EV Outlook 2024 reports that global electric-vehicle (EV) sales approached 14 million units in 2023, a 35% increase over the previous year, lifting the EV share of new-vehicle sales to 18% (IEA, 2024). Driven by this trajectory, charging infrastructure now determines the accessibility, reliability, and emission-substitution benefits that electric mobility delivers to urban residents. Major economies have responded at scale: the European Union's Alternative Fuels Infrastructure Regulation mandates a fast-charging station every 60 km along the trans-European transport network from 2025 (European Parliament and Council, 2023), the U.S. Bipartisan Infrastructure Law allocates 7.5 billion USD to a national charging network (U.S. Department of Transportation, 2022), and China's public charging posts exceeded 4.7 million by the end of 2025, forming the world's largest public network. Yet expansion in deployment scale does not translate automatically into reliable service. As urban climate stress intensifies under extreme heat, heavy precipitation, and high humidity, facility availability risks becoming the weak link of the low-carbon mobility system, and charging faults increasingly constrain operational efficiency, service equity, electrical safety, and the continuity of sustainable urban travel.

The consequences of a charging fault extend well beyond the downtime of a single device. A service interruption immediately causes order losses and station congestion, and a persistent fault can further trigger equipment damage, electrical-safety hazards, and erosion of user trust, discouraging the very adoption of electric mobility on which sustainable cities depend. The deeper problem is that maintenance in the industry remains largely reactive, acting only after a fault has occurred. For an operator managing hundreds to thousands of charging posts, manual inspection cannot cover the complexity of real operation, and reactive repair cannot support refined tasks such as spare-part pre-positioning, staff scheduling, and maintenance-window optimization. If historical transaction logs, meteorological conditions, and equipment attributes up to the current day could instead be used to predict fault risk reliably over multiple future days, maintenance teams could pre-allocate resources and shift the management logic from after-the-fact remediation to preventive intervention. This shift concerns not only the health of individual devices but, more directly, the sustained availability of the charging network, the climate resilience of urban infrastructure, and the stable supply of low-carbon mobility services.

Although research on equipment health management, anomaly detection, and fault prediction has advanced considerably, studies that target long-term, multi-step fault-risk prediction of public charging facilities under climate disturbance, together with the operational response it enables, remain limited in three respects. First, the structural differences among multi-source heterogeneous features are seldom addressed. Charging-fault risk is jointly driven by climate conditions, equipment attributes, user behavior, and fault history, yet these signals differ fundamentally in time scale, statistical distribution, and physical meaning: climate variables act as external disturbances at the daily scale, equipment aging accumulates over months to years, and user behavior and fault history carry short-term fluctuation and autocorrelation. When signals of such different natures are concatenated and learned by a single model, strong predictors such as the one-day-lagged fault count dominate training, whereas weak but physically meaningful slow-varying drivers such as cumulative climate stress are systematically suppressed. This signal-competition problem is rarely treated in a targeted way. Second, existing methods balance multi-scale temporal dependence and long-horizon prediction poorly. Short-term risk is shaped mainly by recent anomalies, whereas medium- and long-term risk is tied to aging, seasonal change, and evolving usage patterns, so that no single time scale describes both. Most studies stop at single-step prediction or short-term anomaly detection and cannot output a daily risk sequence spanning the coming one to four weeks, leaving a clear gap from the resource pre-allocation that real operation requires. Capturing short-term anomalous pulses and long-term evolution within one framework, while remaining robust across horizons, thus remains an open methodological challenge. Third, the step from model output to actionable operation is still missing. A high-accuracy model can output risk probabilities but cannot explain which signals drive an alert, still less quantify how much a given environmental or physical factor changes the fault probability of a specific charging post. Without such explanation and causal identification, predictions are difficult to convert into differentiated, climate-adaptive maintenance. SHAP-based attribution partly reveals which variables matter, but a gap persists between correlational attribution and quantified causal effect, which prevents predictions from offering quantitative support for urban adaptation levers such as summer derating, anti-moisture sealing inspection, spare-part pre-positioning, and differentiated peak-shaving pricing.

To address these problems, this study proposes a multi-step fault risk prediction framework for resilient charging-facility operation, termed FGDSE (Feature-Governed Dynamic Stacking Ensemble). The core idea is as follows. Raw transaction data are first aggregated into post-day time-series samples, and four feature families—physical-environmental, user-behavioral, operational-contextual, and fault-historical—are constructed under strict temporal constraints. Guided by feature partitioning, shallow expert models whose inductive biases match each feature family are then assigned, and two deep experts—a GRU and an enhanced multi-scale temporal convolutional network—are introduced to capture short-term anomalous pulses and long-term evolution simultaneously. At the fusion layer, a horizon-wise dynamic gating mechanism learns adaptive expert weights for each prediction horizon, producing a daily fault risk sequence for the next 1 to 30 days. At the interpretation layer, SHAP and an X-learner extend the analysis from correlation-driven risk scoring to the identification of heterogeneous treatment effects with interventional meaning. The main contributions of this study are as follows:

- This study targets long-term, multi-step fault risk prediction of charging facilities under climate disturbance. It extends the conventional single-day binary classification task into the prediction of a multi-day risk sequence, aligning the model objective with real operational scenarios such as spare-part dispatch, maintenance scheduling, and monthly resource allocation, and addressing the prevailing focus on single-step prediction or after-the-fact diagnosis in existing work.

- The proposed framework is a feature-governed dynamic stacking ensemble. The framework splits the heterogeneous feature space into subdomains by domain semantics and models each with a matched expert, thereby relieving the systematic suppression of weak signals by strong ones. Through an enhanced multi-scale temporal network and a horizon-wise dynamic gating mechanism, it strengthens the characterization of long-horizon temporal dependence and improves cross-horizon adaptability.
- Beyond prediction, this study integrates SHAP attribution and X-learner causal inference into a unified pipeline. It not only identifies the key factors affecting fault risk but also quantifies the heterogeneous treatment effects of different charging posts under environmental and physical disturbances, providing interpretable and actionable evidence for differentiated maintenance.
- Empirically, this study conducts systematic experiments on real operational data covering 13 charging stations over a 25-month period from December 2019 to December 2021. The results show that the proposed framework holds a stable advantage in macro-averaged performance, long-horizon AUC, and recall, and reveal a regular shift in dominance—from fault-history and behavioral signals in the short term to climate and contextual signals at long horizons—as the prediction horizon extends.

The remainder of this paper is organized as follows. Section 2 reviews related work and summarizes the main limitations of existing methods. Section 3 introduces the data source, time-safe feature engineering, and descriptive analysis. Section 4 presents the problem definition and the detailed design of the FGDSE framework. Section 5 reports the experimental results, including performance comparison, ablation analysis, SHAP interpretation, and X-learner causal inference. Section 6 concludes the paper and discusses future research directions.

## 2. Literature review

Around charging-facility fault identification, condition warning, and operational decision support, existing research has largely developed along three paths. The first starts from reliability engineering, mechanistic analysis, and statistical modeling, emphasizing interpretable characterization of system structure and failure processes. The second turns to machine learning and deep learning on operational data, aiming to improve fault identification and time-series prediction accuracy. The third introduces model interpretability and causal identification into charging-infrastructure research. Research targeting long-term, multi-step fault risk prediction of charging posts—and its connection to intervention decisions—remains relatively scarce. On this basis, this section reviews and comments on the related literature along the three dimensions above. Table 1 summarizes how representative studies compare with this study in research focus, method, interpretability, and main difference.

**Table 1. Comparison between representative studies and this study.**

| Reference | Research focus | Method | Interpretability | Key findings | Core contribution | Difference from this study |
|---|---|---|---|---|---|---|
| Quddus et al. (2021) | Reliable charging-station network design under uncertainty | Stochastic optimization | High(explicit model) | Reliability is jointly shaped by capacity, demand, and operation | Couples long-term expansion with short-term operation | Planning/steady state, not post-day multi-step risk prediction |
| Konara et al. (2023) | EVSE failure and repair at fast-charging stations | Continuous-time Markov chain | High | Links effective capacity, completion rate, and resources | Analytical reliability model of failure/repair | No data-driven multi-horizon forecasting |
| Motlagh & Li (2026) | Station siting, sizing, and upgrade planning | Systematic review | N/A (review) | Joint optimization and cross-sector integration underexplored | Synthesizes planning literature | Review; offers no predictive method |
| Kuang et al. (2024) | Effect of electricity price on urban EV charging behavior | Empirical modelling with impulse-response analysis | Medium | Price signals drive marked spatial-temporal heterogeneity in charging demand | Links price response to urban charging-demand patterns | City-level demand behavior, not post-day fault-risk prediction |
| Heo & Chang (2024) | City-scale fast-charging station siting | Advantage actor-critic deep reinforcement learning + geospatial analysis | Medium | Traffic flow and solar potential guide optimal siting; grid balance constrains rewards | Links demand and grid constraints to infrastructure-level siting decisions | Targets siting and planning, not equipment-fault probability |
| Hu et al. (2025) | EVCS growth patterns under EV-adoption trajectories | Allocation model + curve fitting + statistical analysis | Medium | Inflection, convergence, and saturation points vary by city type and traffic structure | Links infrastructure growth stages to sustainable urban development | Long-term planning benchmark, not post-day fault-risk prediction |
| Chen et al. (2024) | Open-circuit fault diagnosis of DC posts | Improved S-transform + LightGBM | Medium (feature-level) | 97.04% mean accuracy, noise-robust | Accurate signal-level diagnosis | Signal level, not operational multi-source risk |
| Feng et al. (2025) | Rectifier open-circuit fault diagnosis | SAE-CNN-ISSA-BiLSTM | Low(deep model) | 99.70% accuracy under strong noise | High-accuracy component diagnosis | Simulated faults, not daily operation |
| Yuan et al. (2026) | DC charger fault diagnosis | Streaming clustering + improved RF | Medium | 97.24% accuracy, fast inference | Lightweight online diagnosis | Classification, not multi-step risk |
| Yang et al. (2024) | Abnormal charging-behavior identification | Feature engineering + random forest | Medium(features) | 0.89 accuracy on AC-side current | Behavior-based anomaly detection | Single-step detection, not horizon sequence |
| Hong et al. (2026) | Multi-step battery-temperature prediction | Vehicle grouping + LSTM | Low | Peak-temperature error 1.45%→1.06% | Multi-step deep prediction on post data | Targets battery temperature, not faults |
| Wang et al. (2026) | Public EV charging usage patterns | SOM temporal clustering + Bayesian spatial multinomial regression | High (spatial model) | Five temporal clusters; functional similarity need not follow geographic proximity | Integrates temporal usage patterns with spatial explanatory modeling | Usage-pattern analysis, not equipment-fault forecasting |
| Hamim et al. (2026) | Real-world EVCS usage and spatial heterogeneity | Traditional models + geographically weighted random forest | High (spatially-aware ML) | Provider, corridor, socioeconomic, and capacity factors drive heterogeneous usage | Uses data fusion to explain station-level utilization across space | Station usage, not post-day fault risk |
| Qiao et al. (2026) | Spatial heterogeneity of EVCS utilization factors | XGBoost + SHAP | High (SHAP) | Fault rate and accessibility are common drivers; brand, agglomeration, and pricing effects vary by area | Explains area-specific utilization mechanisms | Explains utilization, but does not estimate future fault probability or CATE |
| Lu et al. (2025) | Charging-network resilience under rainfall | Resilience curve + XGBoost + SHAP | High (SHAP) | Station count, land use, rainfall drive resilience | Interpretable resilience analysis | Correlational, no causal-effect quantification |
| Zhu et al. (2025) | Causal analysis in transport systems | Causal AI framework | High (causal) | Causal effects identified in transport | Methodological basis for causal AI | Not applied to charging-fault analysis |

| This study | Post-day multi-step fault risk + intervention | FGDSE + SHAP + X-learner | High (SHAP + causal) | Stable long-horizon recall; ~30% heat-sensitive posts | Unifies multi-source forecasting and causal effects | — |
|---|---|---|---|---|---|---|

## 2.1 Charging-facility research based on reliability analysis and statistical models

Early work on charging-infrastructure reliability started mainly from system structure and uncertainty modeling. Quddus et al. (2021) incorporated long-term expansion decisions and short-term operational decisions into reliable station-network design under uncertainty, noting that station reliability is determined not by equipment failure alone but jointly by capacity configuration, demand disturbance, and operational constraints. Konara et al. (2023) explicitly considered the failure and repair processes of EVSE, analyzing the relationship among effective capacity, charging completion rate, and resource coordination at fast-charging stations within a continuous-time Markov-chain framework. At the planning level, related studies have addressed siting-sizing co-optimization, congestion-robust planning, incremental intelligent expansion, and growth-stage benchmarking under rising EV adoption (Motlagh and Li, 2026; Li, Y. et al., 2025; Bai et al., 2025; Golsefidi et al., 2023; Hu et al., 2025). These studies jointly establish an important insight: charging-facility faults are not isolated events, and their consequences propagate to the system level along a chain of capacity contraction, declining service availability, and rising operational cost.

Recent Sustainable Cities and Society studies further show that charging-facility performance is shaped by spatial, economic, and energy-supply contexts. Kuang et al. (2024) examined how electricity price affects urban EV charging behavior in Shenzhen and revealed pronounced spatiotemporal heterogeneity and demand spillover, whereas Heo and Chang (2024) coupled advantage actor-critic deep reinforcement learning with geospatial analysis to optimize city-scale fast-charging station siting under traffic-demand, solar-potential, and grid constraints. These studies show that charging-facility performance is shaped by demand uncertainty and spatial heterogeneity, yet their analytical unit remains the city, network, or station rather than the individual charging post. For preventive maintenance, the key missing link is therefore not only where stations should be built or how demand shifts, but how post-level fault probability evolves under climate and operational disturbance.

The strength of these methods lies in transparent structure and clear parameter meaning, which suits questions such as why a system is vulnerable and how a station network should expand. For the research task of this paper, however, two limitations exist. First, they are better at reliability assessment and expansion planning and cannot directly absorb the high-dimensional heterogeneous information jointly formed by transaction logs, weather disturbance, user behavior, and historical faults. Second, most existing work remains at station-level or system-level steady-state analysis and has not yet formed a prediction framework that operates at the single-post-day granularity and outputs a multi-day future risk sequence. This limitation is the direct motivation for introducing data-driven prediction methods in this paper.

## 2.2 Charging-facility fault diagnosis based on data-driven and deep learning methods

As monitoring and transaction-system data accumulate, the research focus has gradually shifted to data-driven methods. At the fault-identification level, Wang et al. (2023) compared combinations of several word-embedding models and classifiers on real work-order text, showing that operational text records can also support fault-category identification. At the component-level fault-signal diagnosis level, researchers have successively applied deep belief networks (Gao and Lin, 2021), an improved S-transform with LightGBM (Chen et al., 2024), and SAE-CNN-ISSA-BiLSTM (Feng et al., 2025) to the fault diagnosis of DC charging posts, raising diagnostic accuracy from 97.04% to 99.70%. In the battery-system domain, multi-scale feature fusion (Wang, P. et al., 2024) and Transformer-enhanced multimodal sensing (Santhi Mary Antony et al., 2026) have also been used for early fault warning and precise localization. This body of work shows that once fault diagnosis is cast as a structured-feature classification task, machine learning and deep learning offer clear advantages in identification accuracy.

A closer data-driven line examines charging-side behavior, battery states, and charging usage. Yang et al. (2024) used AC-side current features to detect abnormal charging behavior, showing that low-level charging signals can reveal operating anomalies before they become explicit faults. Hong et al. (2025, 2026) further used real-world charging-pile data and vehicle-identification grouping to predict battery voltage and temperature over multiple future steps, demonstrating the value of sequence learning on vehicle-pile interaction data. At the station or network level, charging load and demand have been forecast using LSTM-based multi-feature fusion and federated learning (Aduama et al., 2023; Yin & Ji, 2025), while Wang et al. (2026) integrated temporal clustering with Bayesian spatial multinomial regression to identify public-charging usage patterns in the Bay Area, and Hamim et al. (2026) used real-world data and geographically weighted random forest regression to reveal spatial heterogeneity in EVCS utilization. Recent reviews also emphasize AI-enabled charging control and grid-friendly operation (Li, J. et al., 2026; Ranjbar et al., 2026). These studies confirm that charging data support intelligent prediction and control, but most target behavior, demand, or battery states rather than daily post-level fault-risk forecasting under climate disturbance and the preventive maintenance decisions needed here.

Despite these fruitful results, the existing data-driven literature still differs from the research objective of this paper in three respects. First, the vast majority of work remains at post-fault classification or single-session anomaly detection, rather than long-term multi-step prediction of how fault probability evolves over future days. Second, many studies focus on component-, rectifier-, or battery-level signals, whose prediction targets differ from the daily operational fault risk of public charging posts considered here. Third, existing models organize multi-source heterogeneous features rather coarsely, often concatenating equipment attributes, weather variables, behavioral intensity, and historical faults directly into a single model, and rarely handle explicitly the differences of these signals in time scale and statistical distribution. The feature-governed strategy proposed in this paper is a direct response to the third limitation.

### 2.3 Research on model interpretability, resilience, and causal inference

Beyond prediction accuracy, researchers have begun to ask what a model actually captures from data. Lu et al. (2025) combined resilience curves with XGBoost and SHAP under extreme-rainfall scenarios, identifying station count, the proportion of commercial land, and rainfall conditions as the main determinants of the resilience pattern. Qiao et al. (2026) further used XGBoost and SHAP to explain the spatially heterogeneous determinants of station utilization across urban, suburban, and county areas, showing that fault rate and accessibility are common drivers while brand, agglomeration, and pricing effects vary by spatial context. SHAP has also been used for interpretive analysis of station-selection behavior (Ullah et al., 2023), EV purchase intention (Chen et al., 2026), and charging spatiotemporal patterns (Li, J. et al., 2025). In resilience assessment, Bayesian networks (Guo et al., 2026) and data-driven simulation (Wang, F. et al., 2024) have been used to analyze the vulnerability of charging infrastructure under natural disasters and cyberattacks. This work reflects a common trend: charging-infrastructure research is moving from comparing prediction accuracy alone toward explaining mechanism differences as well.

However, interpretability-oriented studies remain mostly explanatory rather than interventional. SHAP rankings and spatial heterogeneity analyses can clarify which variables are associated with utilization or resilience, but they do not estimate how an environmental or physical exposure would change the future fault probability of a given charging post under a counterfactual intervention. Recent causal-AI work in low-carbon transport shows the value of moving from prediction to policy-effect estimation (Zhu et al., 2025), but comparable causal estimation for charging-post fault risk is still rare. This gap motivates the combination of SHAP attribution and X-learner causal inference in this study.

In summary, existing research has made progress in system reliability assessment, spatially heterogeneous charging demand and utilization analysis, fault diagnosis, short-term time-series prediction, and interpretable learning. For the daily-granularity operational scenario of charging posts, however, a unified analytical framework is still missing: one that can simultaneously absorb multi-source heterogeneous signals, output long-term multi-step fault risk, and quantify the heterogeneous causal effects of environmental and physical factors. This study proposes the FGDSE framework to connect the transport problem of resilient charging-infrastructure operation with the engineering problem of device-level preventive maintenance.

## 3. Data collection and analysis

### 3.1 Data source

The charging-facility operational data used in this study come from 13 charging stations under the Tongxiang Power Supply Company of State Grid Zhejiang Electric Power Co., Ltd.; the raw dataset was publicly released by Zhang et al. (2025) in 2025. The data span December 2019 to December 2021 and cover 441,077 charging-session records. Each record contains the charging-post ID, user ID, order creation time, charging energy, electricity cost, service charge, charging start and end times, termination reason, and meteorological variables such as daily mean temperature, relative humidity, and precipitation. Data acquisition followed the IEC 60870-5-104 protocol with dual Bluetooth and 4G communication, and transmission used SSL/TLS encryption, ensuring data integrity and security.

The raw data fully cover the period from the outbreak of COVID-19 to normalized recovery. According to the dataset documentation, the first quarter of 2020 was affected by travel controls and restricted operation, and its anomaly-flag rate (0.2504) was markedly higher than that of the recovery period (0.1902). Considering the sustained disturbance to travel from localized epidemic-control measures during 2020–2021, and to preserve the integrity of the time series while objectively testing the model's predictive robustness under non-stationary conditions, this study retains the full dataset for training and evaluation without manual exclusion. This treatment aims to verify the generalization ability of the FGDSE framework under exogenous shocks. The final sample spans December 2019 to December 2021 and contains 57,237 post-day observation records for 92 charging posts. Table 2 lists the main fields of the dataset.

**Table 2. Main fields of the dataset.**

| Field | Description |
|---|---|
| User ID | Anonymized unique identifier of the user |
| Charging Post ID | Unique identifier of the charging post |
| Order Creation Time | Order creation timestamp(to the second) |
| Transaction Power (kWh) | Charging energy of the session |
| Electricity Cost/Service Charge/Transaction Amount | Electricity cost, service charge, and total transaction amount |
| Start Time/End Time | Charging start and end times |
| Is_nonsystem_fault | Fault flag of the main or control circuit(binary variable) |
| Is_user_stop/Is_full_charge/Is_abnormal | Other termination reasons and a composite anomaly flag |
| Temperature(°C)/RelativeHumidity(%)/Precipitation(mm) | Daily meteorological variables |
| Location Information | Location category of the charging station |

### *3.2 Data processing*

#### *3.2.1 Data cleaning*

The dataset had undergone preliminary cleaning at release, including desensitization of user-identification information, removal of incomplete records (about 0.022%), and deletion of sessions with abnormal charging duration, energy exceeding the rated capacity, or clearly anomalous payment amounts. To assess quantitatively whether cleaning introduced distribution shift, this study used a two-sample Kolmogorov–Smirnov (K-S) test to compare the empirical cumulative distributions before and after cleaning. The K-S test is a non-parametric test of distributional consistency; its null hypothesis is that the samples before and after cleaning come from the same distribution, and the statistic D is the maximum vertical distance between the two empirical cumulative distribution functions, with a smaller D indicating a weaker difference. The K-S results on the raw dataset show extremely small distribution differences in three key indicators—charging energy (D=0.00006), cost (D=0.00006), and charging start time (D=0.00004)—indicating that cleaning mainly removed obvious anomalous records without substantially altering the empirical distributions of the key variables.

On this basis, this study completed a key task: spatiotemporal alignment of transaction data and meteorological data. Charging-transaction data are recorded per session, whereas meteorological data are daily averages for the administrative district of the station, so the two differ in temporal granularity. The alignment strategy uses the charging-post ID and date as a composite key and merges the daily meteorological data into the transaction data through a left join. The left join preserves all transaction records and avoids discarding high-value samples due to local gaps in the meteorological data. For the very few dates with missing meteorological measurements, linear interpolation between adjacent dates was applied, with the continuous gap not exceeding three days. After alignment, the transaction data were aggregated by charging-post ID and date to produce a daily statistical table at the post-day granularity, with aggregated indicators including daily total charging energy, daily total transaction count, daily mean charging power, number of non-system faults (i.e., main- or control-circuit faults), number of user interruptions, and date attributes. The fault label is defined as follows: if a charging post experiences at least one non-system fault on day t, then $\mathrm{y_t} = 1$; otherwise $\mathrm{y_t} = 0$.

### *3.2.2 Feature engineering*

The core constraint of feature construction is that features used to predict the state on day $\mathrm{t}$ may use only information up to day $\mathrm{t} - 1$ and earlier. This is not merely a technical requirement but the logical premise of causal inference. All features involving historical statistics, such as rolling means and cumulative sums, are implemented with a one-step lag. Continuous features are standardized with Welford's algorithm independently on each charging post's time series, and the standardization parameters are computed only from data before the current time, eliminating future-information leakage at the source.

Based on domain knowledge and the physical failure mechanisms of charging facilities, this study divides the features into four mutually exclusive feature families.

**Physical feature family.** This family describes the external effect of climate conditions on charging equipment, including temperature, humidity, precipitation, and their statistics. The dew-point depression ΔD is introduced as a proxy indicator of surface-condensation risk on electrical equipment. The dew-point temperature $T_d$ is computed with the Magnus–Tetens approximation:

$$T_d = \frac{c \cdot \gamma(T, RH)}{b - \gamma(T, RH)}, \gamma(T, RH) = \frac{b \cdot T}{c + T} + \ln\left(\frac{RH}{100}\right) \tag{1}$$

where $b = 17.67$ and $c = 243.5$ are empirical constants, $T$ is the Celsius temperature, and $RH$ is the relative humidity (%). The dew-point depression is defined as $\Delta D = T - T_d$. When $\Delta D < 3°C$, the surface-condensation risk is high, and moisture ingress into the insulating material may lead to electrical failure. In addition, a Meteorological Weighted Index (MWI) is constructed to characterize composite climate stress:

$$MWI_t = \frac{w_T \cdot \tilde{T}_t + \tilde{H}_t + w_R \cdot \tilde{R}_t}{3} \tag{2}$$

where $\tilde{T}_t$ , $\tilde{H}$ , and $\tilde{R}_t$ are the historical-quantile normalized values of temperature, humidity, and precipitation, respectively; $w_T$ takes 1.5 for extreme high temperature (>35°C) and extreme low temperature (<5°C) and 1.0 otherwise; and $w_R$ takes 2.0 on rainy days and 1.0 otherwise. The above physical variables are further summarized by mean, maximum, and standard deviation within sliding windows of 3, 7, 14, and 30 days, forming multi-scale climate-exposure features that characterize the cumulative effect of sustained climate stress.

**User-behavioral feature family.** This family characterizes the usage pattern of a charging post and reflects the dynamic effect of user behavior on equipment health. Core features include the load intensity $L_t = E_t/N_t$ (average energy per session), the share of transactions in the peak period (19:00–21:00), the user repeat-use rate (the ratio of total visits to unique users), and the zero-power charging share (the proportion of invalid sessions caused by abnormal connector contact). These variables are likewise rolled over multiple time windows. This family captures the dimension of how users use the equipment, whose physical meaning corresponds to the cumulative effect of mechanical wear and thermal fatigue.

**Operational-contextual feature family.** This family contains the static attributes and slowly varying operational background of a charging post. Equipment age is measured by the number of operating days since the post first appears in the dataset and is discretized by tertiles into new, middle-aged, and old equipment, corresponding to the physical stages of material aging. Other features include the regional historical mean fault rate, the ratio of the post's historical fault rate to the regional mean, and the cross-cumulative mean of region and day-of-week.

**Fault-historical feature family.** This family uses as its core features the one-day-lagged fault count and its rolling mean, maximum, and standard deviation within 3-, 7-, 14-, and 30-day windows. The cumulative mean fault rate of a charging post reflects its baseline fault tendency. To capture fault-acceleration trends, the Trend Stress Index (TSI) is constructed. Its computation first fits a linear trend coefficient $\beta_t$ within a sliding window $W$:

$$\beta_t = \text{LinearTrend}(F_{t-W+1}, \dots, F_t) \tag{3}$$

$$TSI_t = \frac{\beta_t}{\overline{F}_t + \epsilon}, \overline{F}_t = \frac{1}{W}\sum_{i=0}^{W-1} F_{t-i} \tag{4}$$

where $\epsilon$ is a small constant preventing a zero denominator, and $W = 7$ days. $TSI > 0$ indicates that the fault frequency is accelerating—a leading risk signal that a static mean cannot reflect. In addition, several interaction features with clear physical or business meaning are constructed; for example, inter_temp_fault = temp_max × fault_count_lag1 characterizes the coupling between high temperature and recent faults, and inter_phy_current_stress = mwi × fault_count_lag1 captures the synergy between composite meteorological stress (MWI) and recent fault history.

### *3.2.3 Strict causal constraint and online standardization*

Continuous features are standardized with Welford's online algorithm independently on each charging post's time series:

$$\mu_{s,t} = \mu_{s,t-1} + \frac{x_{s,t} - \mu_{s,t-1}}{n_t} \tag{5}$$

$$\sigma^2_{s,t} = \sigma^2_{s,t-1} + \left(x_{s,t} - \mu_{s,t-1}\right)\left(x_{s,t} - \mu_{s,t}\right) \tag{6}$$

where $n_t$ is the number of samples up to time $t$. The standardization parameters $(\mu_{s,t-1}, \sigma_{s,t-1})$ used to predict day-$t$ features are computed only from data up to day $t-1$ and earlier, introducing no future information at any point and ensuring the temporal safety of feature engineering.

## 3.3 Analysis of climate-influenced charging behavior and fault characteristics

Before building the model, this section first conducts a descriptive analysis of the relationships among climate variables, charging behavior, and fault risk, providing an empirical basis for the subsequent feature selection and causal inference.

**Nonlinear relationship between temperature and charging demand.** Figure 1 shows the scatter distribution and smoothed trend of daily mean temperature against the network-wide daily total of charging transactions.

The relationship is not monotonic but exhibits a nonlinear, valley-shaped pattern: in the low-temperature range (around −3°C) the daily transaction total stays relatively high (about 900 per day); it falls to a trough of about 430 per day in the mid-temperature range (about 8–9°C); and as temperature rises further it gradually recovers, holding a plateau of about 580–600 per day in the 15–24°C range. In the high-temperature range above 30°C the transaction count does not decline but climbs further to about 880 per day, indicating that the suppressive effect of high temperature on charging demand is limited and far weaker than its aggravating effect on equipment thermal stress. This asymmetric, non-monotonic influence of temperature on the demand side and the equipment side provides empirical support for introducing multi-scale temperature rolling features into the model.

**Electricity pricing also shapes the spatiotemporal distribution of charging behavior.** In Jiaxing, the off-peak price (11:00–13:00 and 22:00–08:00 the next day) is 0.3784 CNY/kWh, only 31% of the peak price (19:00–21:00, 1.2064 CNY/kWh). The data show that charging demand in the off-peak period is relatively concentrated, and users exhibit clear price-responsive behavior. This implies that charging load is not randomly distributed but is markedly regulated by economic incentives. The degree of load concentration is directly related to the fatigue-accumulation rate of charging equipment: dense charging in peak periods may cause local thermal-load superposition and impose higher demands on the heat-dissipation structure. Based on this observation, this study includes the peak-period transaction share in the behavioral feature family to capture the effect of price-driven charging concentration on equipment health.

**Precipitation has complex nonlinear effects on charging behavior and equipment condition.** Figure 2(a) compares per-post daily charging counts under clear and extreme-rain conditions. The sample distribution shows that, although clear-day data are abundant and contain extreme high-frequency records, a large number of low values pull down the overall median; the extreme-rain sample is smaller, but its charging-frequency distribution is more concentrated and its overall center shifts upward markedly, with a statistically marginal effect (p=0.0669).

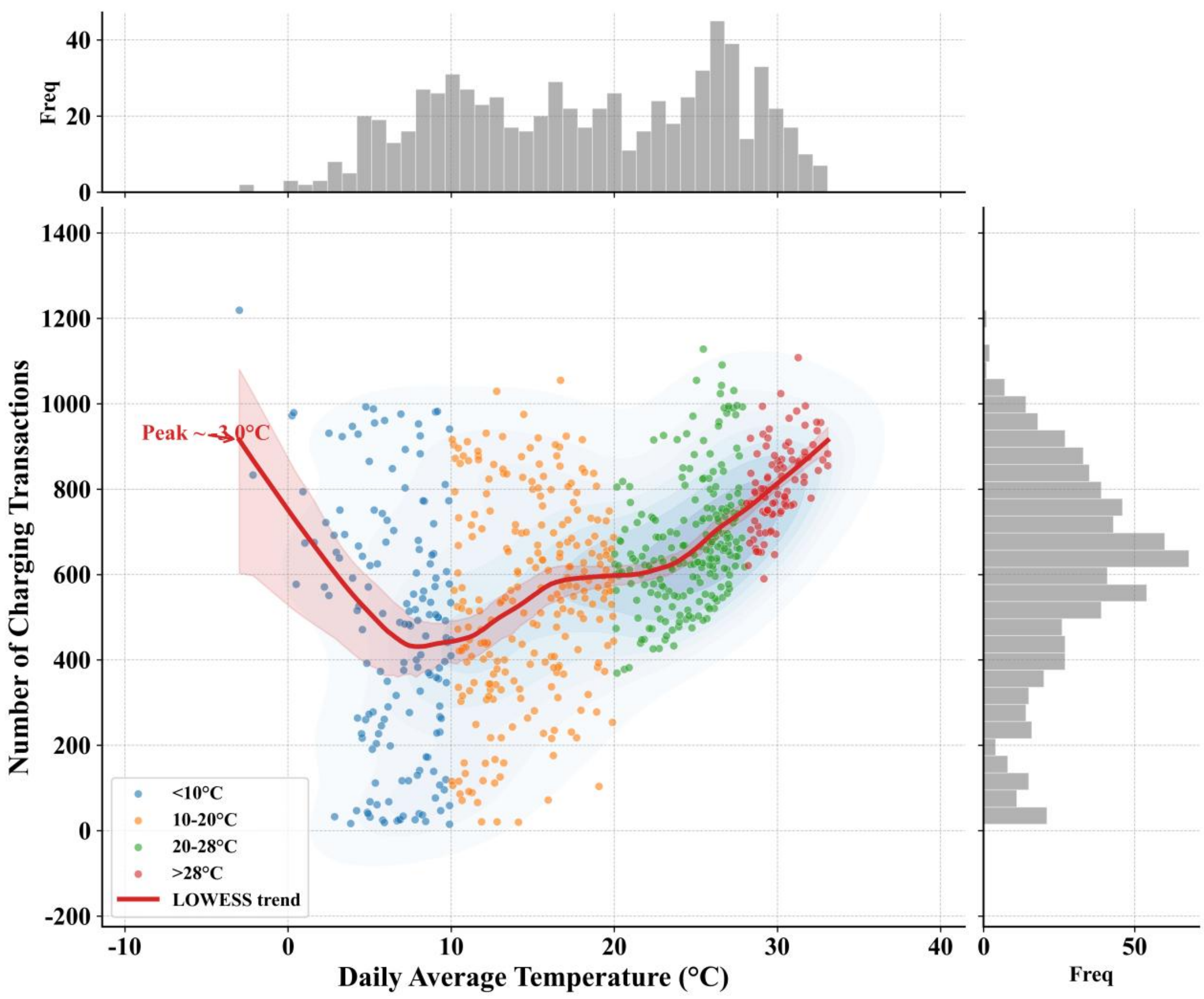

**Figure 1.** Daily mean temperature versus number of charging transactions.

This indicates that adverse weather does not substantially reduce overall charging demand. Heavy rain raises driving energy consumption and forces vehicles still in service to increase daily recharging. This rise in per-vehicle charging frequency offsets the reduction in total vehicles on the road, so charging facilities still bear a high operational load under extreme weather. Figure 2(b) further analyzes changes in zero-power charging events after precipitation: the counts of zero-power events on the rain day and on days 1–3 afterward show no significant difference (Wilcoxon test $p > 0.6$), and the ratio curve stays below 1.0 throughout. This means that no precipitation-induced moisture-lag effect on connectors is observed in this dataset. Even so, this study retains the precipitation rolling features to capture possible long-term cumulative effects, because a weak effect not observable in the short term may become significant over longer time scales.

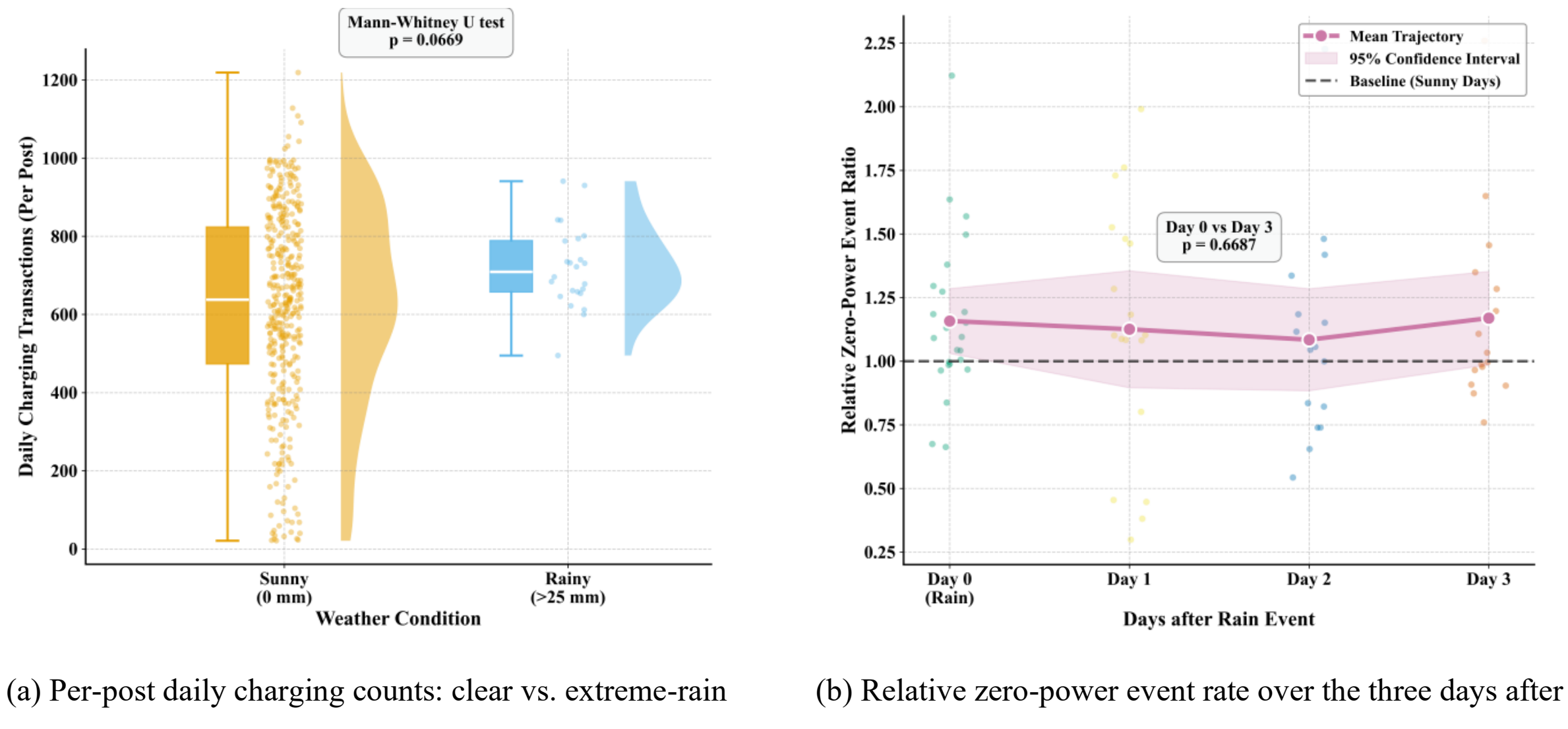

(a) Per-post daily charging counts: clear vs. extreme-rain conditions

(b) Relative zero-power event rate over the three days after precipitation

**Figure 2.** Effect of precipitation on charging behavior and zero-power charging events.

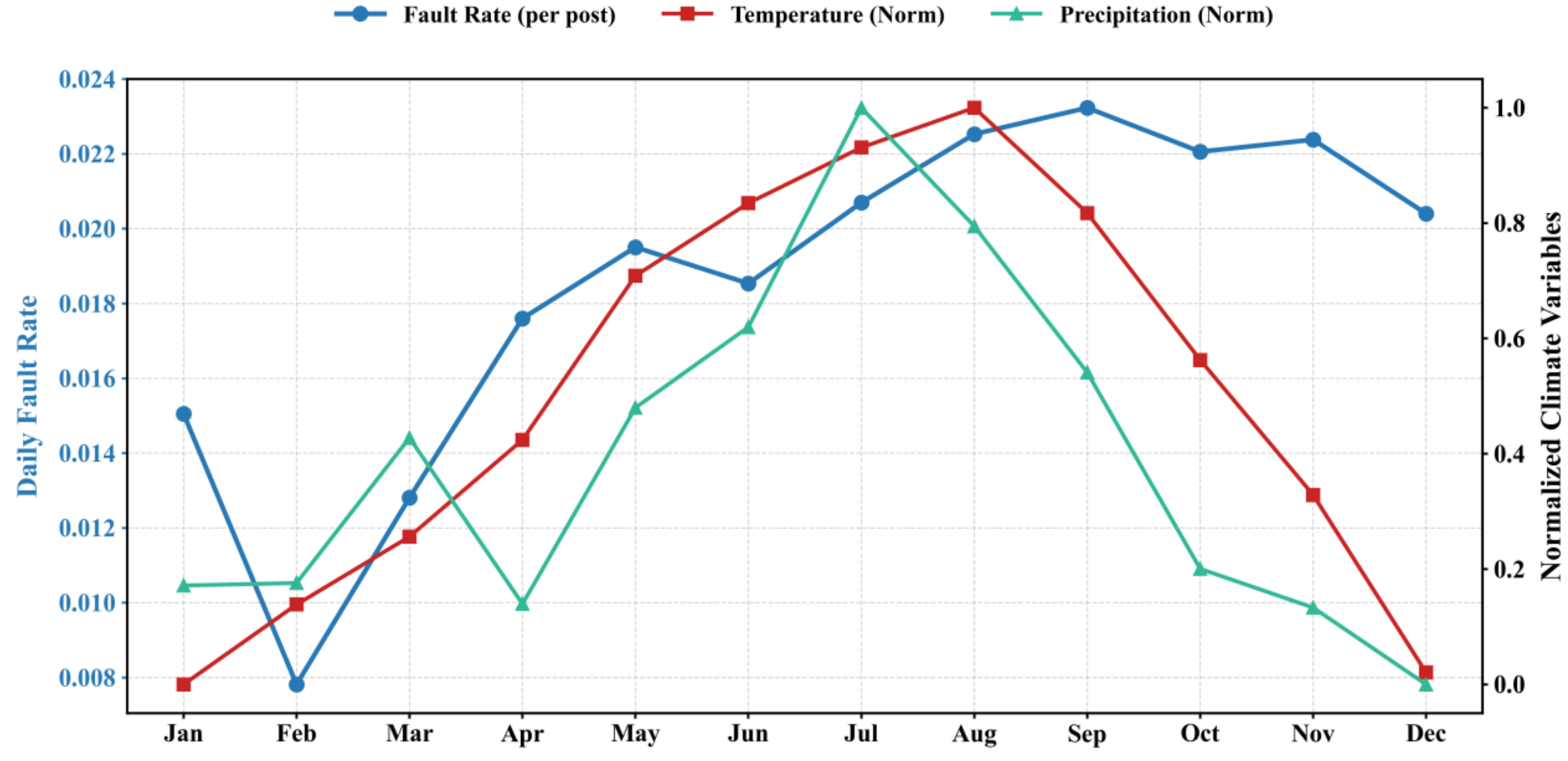

**Figure 3.** Seasonal co-variation of fault rate and climate variables.

**Seasonal co-variation of climate and fault rate.** Figure 3 compares the monthly variation of fault rate, temperature, and precipitation. The fault rate reaches its annual peak in September (about 0.023 per day), showing a marginally significant positive association with temperature ($r = 0.572, p = 0.052$), whereas no significant linear association with precipitation is observed ($r = 0.284, p = 0.371$); yet the asynchronous fluctuation of the two at the monthly scale suggests a potential nonlinear coupling, and this seasonal regularity persists across the two years. From the physical-mechanism perspective, the thermal-fatigue rate of the power-conversion module increases exponentially with junction temperature, the connector interface tends to become unstable under sustained high temperature due to metal creep, and the aging rate of insulating material rises markedly with temperature following the Arrhenius relation.

A point worth deeper discussion is that the relationship between precipitation and fault rate is nonlinear: the surge of precipitation in June is not accompanied by a synchronous rise in fault rate, reflecting the mutual offset between the fault-inducing effect of precipitation (connector moisture ingress) and the fault-suppressing effect (reduced travel and thus reduced mechanical wear). This finding reinforces an important judgment: correlation analysis alone cannot accurately characterize the true effect of climate on faults, which is precisely the core motivation for introducing the X-learner causal-inference framework in this study.

**Individual heterogeneity of fault tendency across charging posts.** Figure 4 displays the large differences in fault rate across charging posts: the highest-risk post reaches a daily mean of 2.815 faults, whereas the most stable post is only 0.062, a range of more than 45 times. This difference is not random noise but the combined result of equipment attributes, usage intensity, and operational background.

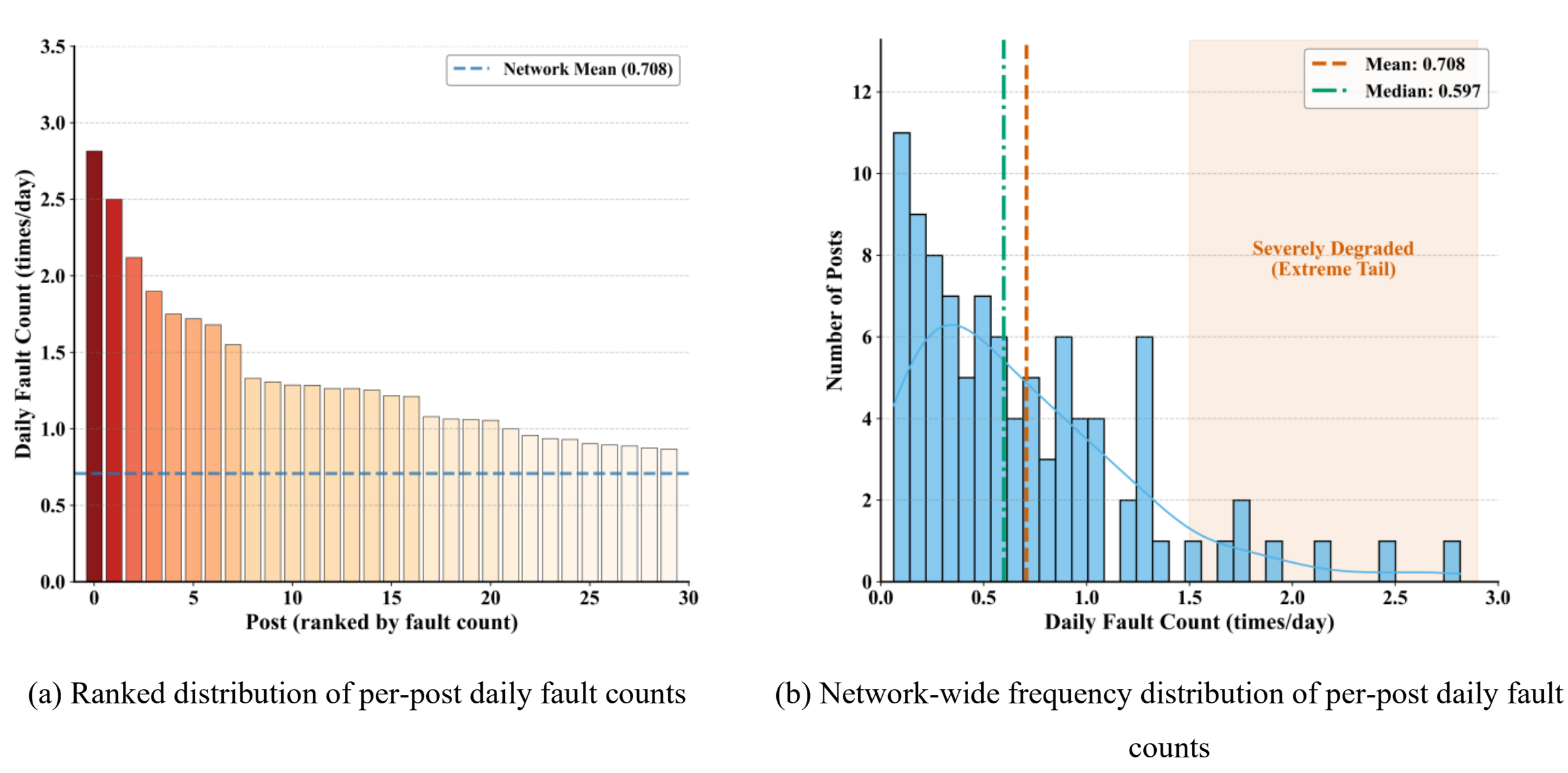

(a) Ranked distribution of per-post daily fault counts

(b) Network-wide frequency distribution of per-post daily fault counts

**Figure 4.** Individual heterogeneity of post-level fault rates.

Figure 5 further reveals the interaction between equipment age and high temperature: equipment is most sensitive to high temperature in its early service (the high-temperature group's daily mean fault count is about 63.9% higher than the low-temperature group), the sensitivity falls markedly after entering steady operation (about +8.5%), and rises again in the aging period (about +22.8%), forming a typical bathtub-shaped sensitivity curve. These patterns provide a prior direction for identifying the heterogeneity of high-temperature treatment effects in the subsequent causal inference.

Taken together, the influence of climate variables on charging-facility fault rate exhibits marked heterogeneity and nonlinearity. Different climate factors act on equipment through different physical pathways, and the effect magnitude varies with the age and usage intensity of the charging post. These patterns form the empirical basis for systematically incorporating climate features into the model and for focusing on climate variables in the subsequent SHAP analysis and X-learner causal inference.

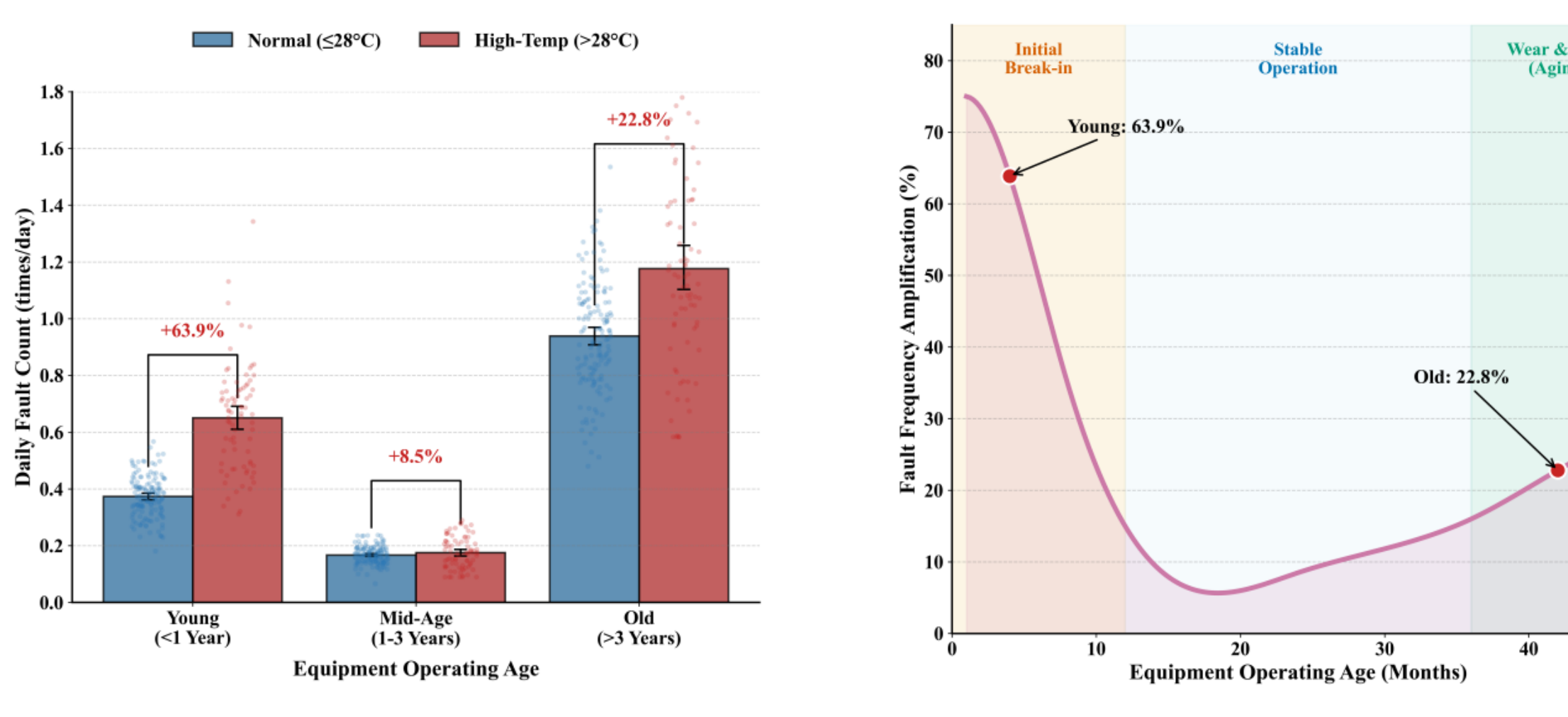

(a) Daily fault counts by equipment-age group (b) Fault-frequency trajectory across the equipment service life

**Figure 5.** Interaction between equipment age and high temperature.

# 4. Methodology

This section presents the complete technical design of the FGDSE framework. We first formalize the problem, then describe the three core stages of the framework in turn: time-safe feature engineering (detailed in Section 3), divide-and-conquer expert prediction, and meta-fusion decision-making. Finally, we introduce the hyperparameter optimization strategy and the causal inference module.

## 4.1 Problem statement

This study focuses on proactive fault warning for electric vehicle supply equipment (EVSE). The core task is to output a daily fault risk sequence over the next H days, based on historical observations up to the current day t.

Let the multivariate input series up to time t be denoted by $X_{1:t} = \{(x_\tau, y_\tau)\}_{\tau=1}^{t}$, where each $x_\tau \in \mathbb{R}^d$ is the d-dimensional feature vector for day τ and each $y_\tau \in \{0,1\}$ is a binary fault indicator, with $y_\tau = 1$ meaning at least one fault occurred on that day. The multi-step prediction problem is then formalized as learning a mapping function $f$:

$$\hat{Y}_{t+1:t+H} = f(X_{t-L+1:t}) \tag{7}$$

$$(\hat{y}_{t+1}, \hat{y}_{t+2}, \dots, \hat{y}_{t+H}) = f(X_{t-L+1:t}) \quad (8)$$

Here the recent feature window has length L, and the predicted value for day t+H is the estimated probability that a fault occurs on that day. This formulation extends fault warning from a single-point judgment into a sequential risk assessment.

## 4.2 Overall framework design

The FGDSE framework adopts a three-stage progressive architecture (Figure 6). Its core idea is to divide and then coordinate.

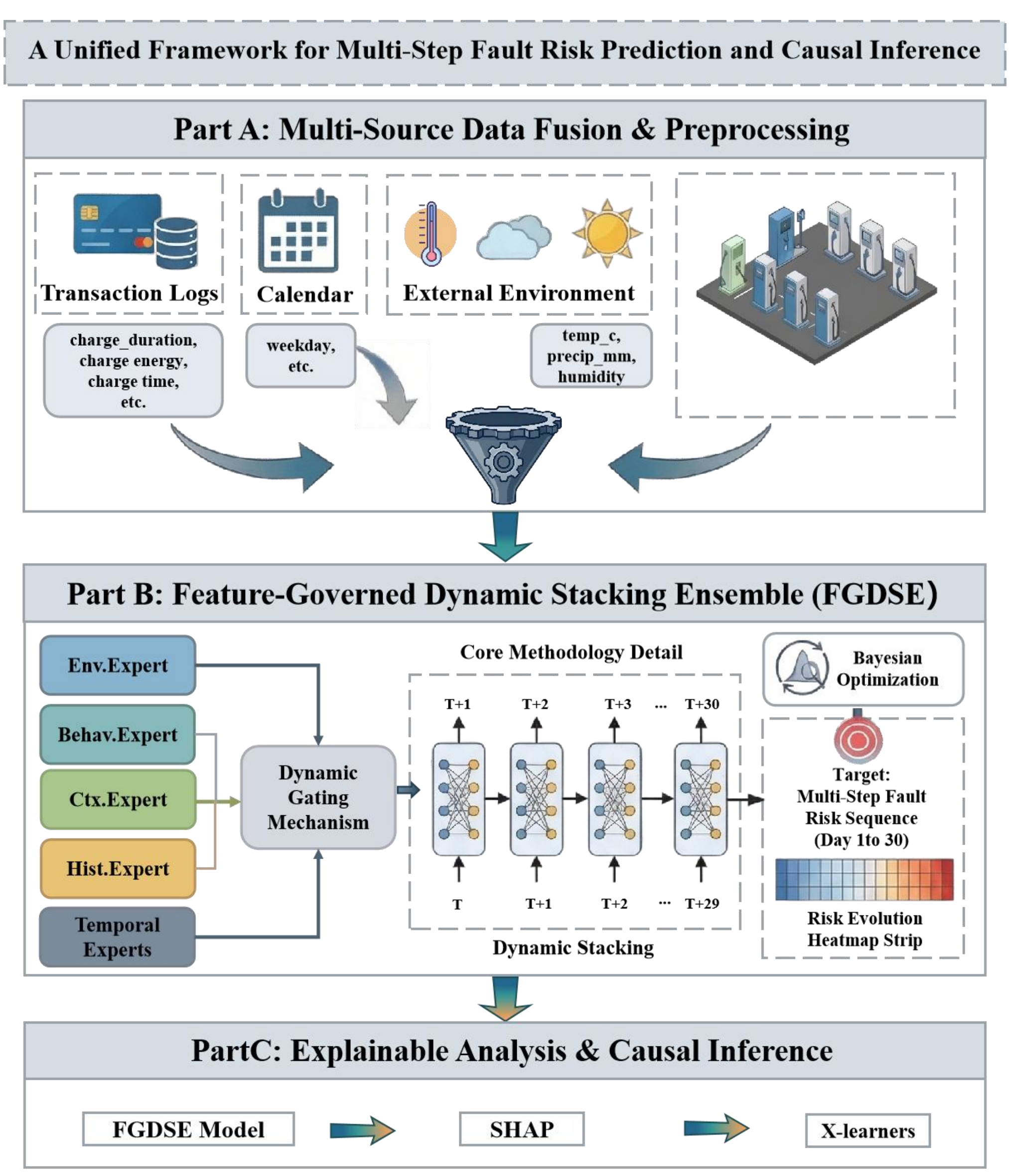

**Figure 6.** The FGDSE multi-step fault prediction framework based on multi-expert fusion.

**Stage 1: Time-safe feature engineering.** Raw transaction logs are aggregated at the daily granularity.

Guided by physical principles and domain knowledge, four feature families are constructed: physical, behavioral, contextual, and historical. All features are built under strict causal constraints, as detailed in Section 3.2.

**Stage 2: Divide-and-conquer expert prediction.** Six experts are designed. Four are shallow experts targeting specific feature families (built on LightGBM, XGBoost, and CatBoost), and two are deep temporal experts (a GRU and an enhanced multi-scale TCN). Each expert independently outputs a vector of fault probabilities over the next H days.

**Stage 3: Meta-fusion decision-making.** A horizon-wise dynamic gating mechanism is proposed. It learns a set of optimal expert weights for each future day, and adaptively balances the reliability of different experts across different prediction horizons.

## 4.3 Partitioned expert predictors

### *4.3.1 Domain experts (shallow models)*

Matching each feature family with a base learner whose inductive bias fits the data is the key step that turns the idea of feature partitioning into practice. Based on the data form and the physical meaning of each family, four differentiated gradient-boosting tree models are selected for the four feature domains.

**Physical expert (LightGBM).** The physical family consists mainly of continuous sensor data. The histogram-based algorithm of LightGBM is both efficient and effective on such data, and its GOSS mechanism accelerates training without loss of accuracy.

**Behavioral expert (XGBoost).** The behavioral family is high-dimensional (about 80 to 100 columns) and contains complex interactions among variables. XGBoost has advantages in handling sparse features and in column-subsampling regularization, which helps capture the nonlinear links between user behavior and fault risk.

**Contextual expert (CatBoost).** The contextual family contains categorical features such as equipment-age bins and regional codes. CatBoost handles categorical features natively, and its symmetric tree structure fits static attribute features more robustly than other GBDT variants.

**Historical expert (XGBoost).** The historical family shows strong autocorrelation and a nonlinear evolution trend. The function-approximation ability of XGBoost captures multi-scale patterns in the fault sequence, from short-term fluctuation to long-term trend. Although the behavioral and historical experts both use XGBoost, their input feature spaces differ in essence; replacing the historical expert with LightGBM or ExtraTrees lowers the AUC by about 0.003 to 0.005 in our experiments.

Each expert takes the data of its feature domain as input and directly outputs an H-dimensional probability vector.

#### *4.3.2 Deep temporal experts*

To make up for the limited ability of shallow models in capturing complex long-term temporal dependence, two deep sequence experts are introduced. Both take the standardized multivariate window as input.

**GRU expert.** The hidden state is updated step by step over time:

$$h_\tau = GRU\left(x_\tau, h_{\tau-1};\ W_{gru}\right),\ \tau = t - L + 1, \dots, t \tag{9}$$

The hidden state at the final step is mapped by a multilayer perceptron into H-step logits, and a sigmoid function then produces the probability output.

**Enhanced multi-scale TCN expert.** Three key enhancements are made to the standard TCN. First, multi-branch parallel convolution: kernels of size three, five, and seven are used in parallel branches to capture short-term pulses, medium-term fluctuation, and long-term trend, respectively. Second, depthwise separable convolution with SE attention: following the idea of combining SE attention with temporal convolution to improve time-series prediction (Wu et al., 2026), a Squeeze-and-Excitation module is added within each branch. Global average pooling compresses the spatiotemporal information into a channel descriptor, and a gating mechanism then produces channel-importance weights:

$$s = \sigma\left(W_2 \delta(W_1 z)\right),\ \ z_c = \frac{1}{L} \sum_{\tau=1}^{L} U_{c,\tau} \tag{10}$$

Here δ is the GeLU activation and σ is the sigmoid function. Third, attention-based temporal pooling: a lightweight attention module adaptively aggregates information from the key time steps:

$$\alpha = Softmax\left(v^T tanh(W_a H + b_a)\right) \tag{11}$$

Here $H \in \mathbb{R}^{C \times L}$ is the feature map after branch fusion, the final representation $r = H\alpha^T$ is then passed to an MLP, which outputs the probability.

## 4.4 Meta-fusion strategy

The six experts produce six groups of H-step probability vectors. The meta-fusion layer aims to learn a combination function that yields a final prediction more accurate than any single expert.

#### *4.4.1 Stacked ensemble*

A separate meta-learner (a lightweight XGBoost) is trained for each horizon k. The predicted probabilities of all experts at horizon k form the feature vector $q_{t,k} = \left(p_{t,k}^1, \dots, p_{t,k}^6\right)$:

$$\hat{y}_{t+k}^{stack} = M_k\left(q_{t,k};\ \Theta_k\right) \tag{12}$$

#### *4.4.2 Horizon-wise dynamic gating*

To pursue higher parameter efficiency and interpretability, a horizon-wise dynamic gating mechanism is proposed. For each horizon k, a fixed expert weight vector is learned, subject to the constraint that the weights sum to one:

$$\hat{y}_{t+k}^{gate} = \sum_{j=1}^{6} w_{k,j} \cdot p_{t,k}^{j} \tag{13}$$

The weight matrix is optimized jointly by minimizing the total binary cross-entropy loss on the validation set:

$$min_W \quad -\frac{1}{N_v H}\sum_{i=1}^{N_v}\sum_{k=1}^{H}\left[y_{i,t+k} \log \hat{y}_{i,t+k}^{gate} + \left(1-y_{i,t+k}\right) \log\left(1-\hat{y}_{i,t+k}^{gate}\right)\right] \tag{14}$$

$$s.t.\ w_k \geq 0,\ \ 1^T w_k = 1,\ \ \forall k \in \{1, \ldots, H\} \tag{15}$$

This is a constrained convex optimization problem, solved efficiently by projected gradient descent. For numerical stability, the direct BCE gradient is used. The learned weight matrix has an intuitive interpretation: each row reveals how much the system trusts each expert when predicting the k-th future day. Figure 7 illustrates the internal structure of the meta-fusion layer and the horizon-wise dynamic gating, including the shape and constraints of the weight matrix and its coupling with the output vector of each expert.

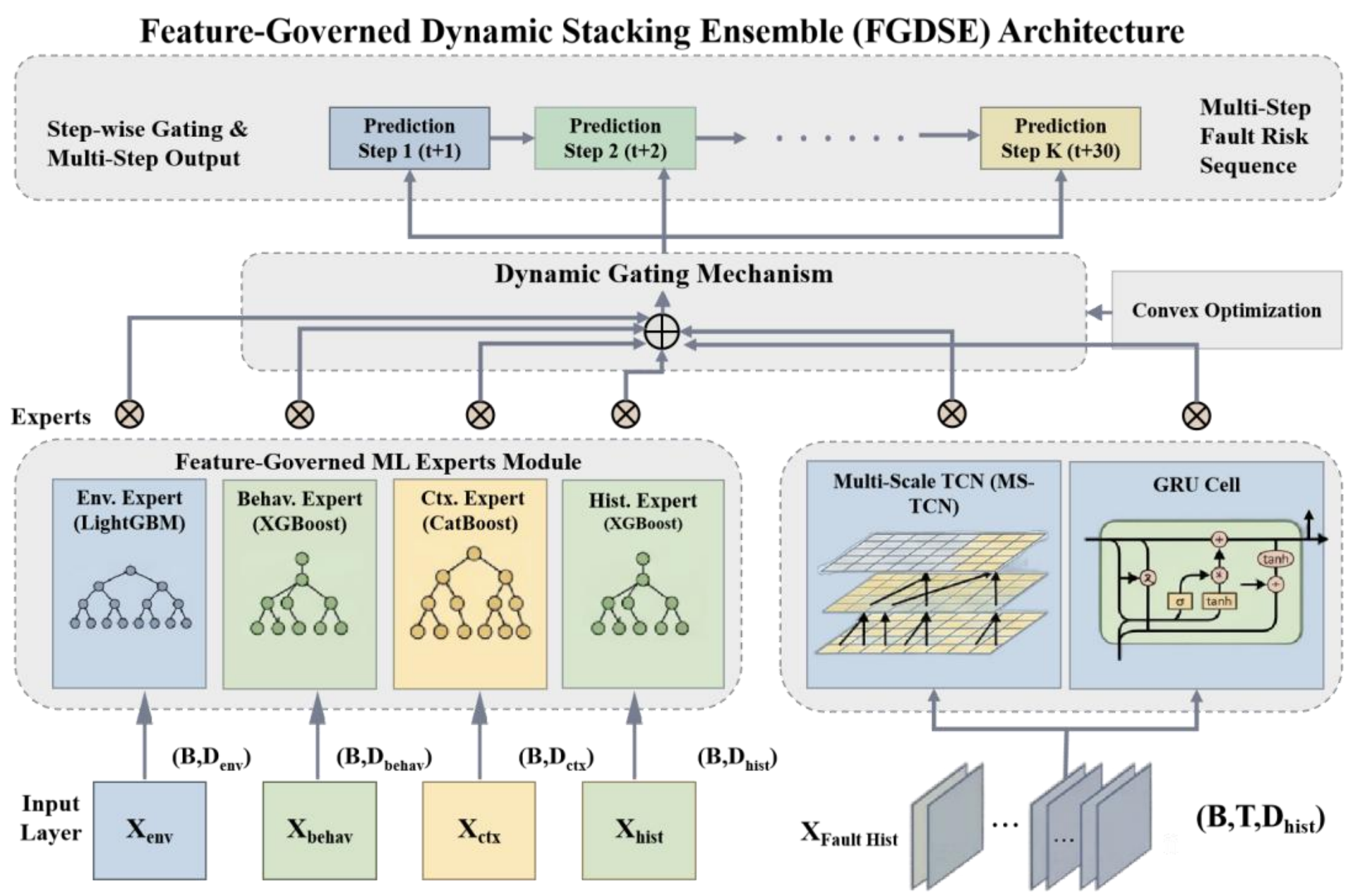

**Figure 7.** Internal structure of the meta-fusion layer and the horizon-wise dynamic gating mechanism.

## 4.5 Hyperparameter optimization

A Bayesian optimization algorithm based on the Tree-structured Parzen Estimator (TPE) is adopted, implemented with the Optuna framework for automated search. Using a target threshold, TPE divides past evaluations into a well-performing group and a poorly-performing group, and builds a probability density function for each:

$$p(x|y) = \begin{cases} l(x) & if\ y < y^* \\ g(x) & if\ y \geq y^* \end{cases} \tag{16}$$

The next configuration to explore is chosen by maximizing the ratio of the two densities. A median pruner is also introduced to stop early any trial that falls below the historical median during training. Bayesian optimization runs only on the training and validation sets; the test set never participates.

## 4.6 Causal inference with the X-learner

An accurate prediction model can flag high-risk charging posts, but it cannot answer the core counterfactual question: if a specific intervention is applied to a certain type of post, how much would the fault probability actually drop? To overcome this limit, the X-learner framework is introduced to estimate, from observational data, the conditional average treatment effect (CATE) of each treatment variable, and to further quantify the average treatment effect (ATE).

### *4.6.1 Treatment variables and key assumptions*

Let the binary treatment variable take values of zero or one. Continuous environmental and operational indicators, such as maximum temperature, precipitation, and peak load, are discretized using a specific quantile of the sample distribution as the threshold. This study uses the 85th percentile for high peak load and the 90th percentile for extreme heat, heavy rain, and high humidity, so as to focus on extreme exposure conditions:

$$T_i = \begin{cases} 1 & if\ V_i > Q_q(V) \\ 0 & otherwise \end{cases} \tag{17}$$

Here the threshold is the q-th quantile of the variable. The CATE is defined as the expected individual treatment effect given the covariates:

$$\tau(x) = E\left[Y^{(1)} - Y^{(0)} \mid X = x\right] \tag{18}$$

The validity of the X-learner relies on the strong ignorability assumption: given the observed covariates, treatment assignment is independent of the potential outcomes. In this study, the covariate set covers four families of features, including equipment age, historical fault rate, rolling climate statistics, and user-behavior patterns, which adjust fairly well for systematic differences between the treatment and control groups. Even so, supply-side responses such as real-time inspection and temporary power derating by operators are not directly included as covariates, and may form unobserved confounding. For this reason, several diagnostics and sensitivity checks are added outside the estimation procedure to assess how robust the conclusions are to potential violations of this assumption.

### *4.6.2 Propensity score estimation and the common-support check*

The propensity score is estimated by logistic regression with class-balancing weights. To check whether the treatment and control groups are comparable in the covariate space, two diagnostics are performed.

**(1) Overlap of the propensity score distributions.** The probability density curves of the two groups are plotted, and the overlap ratio, defined as the probability mass jointly covered by both distributions, is computed. An overlap ratio closer to one indicates a more adequate common support.

**(2) Covariate balance check.** The standardized mean difference (SMD) is used to quantify the between-group difference of each covariate before and after weighting. For a continuous covariate, the SMD is defined as:

$$SMD_j = \frac{X_j^1 - X_j^0}{\sqrt{\frac{s_{j1}^2 + s_{j0}^2}{2}}} \tag{19}$$

Here the means and standard deviations are taken over the treatment group and the control group. Inverse probability weighting (IPW) is applied to adjust for the propensity score, and the maximum SMD before and after weighting is reported. A weighted maximum SMD below 0.1 is generally regarded as an acceptable level of covariate balance.

#### *4.6.3 X-learner estimation procedure*

The X-learner has three steps, and LightGBM is used as the base learner.

**Step 1: Response-function estimation.** The conditional expectation models for the control and treatment groups are fitted separately:

$$\mu_0(x) = E[Y \mid T = 0, X = x],\ \ \mu_1(x) = E[Y \mid T = 1, X = x] \tag{20}$$

**Step 2: Imputed-effect computation.** Cross-prediction is used to construct proxies for the individual treatment effect:

$$\tilde{\tau}_1^i = Y_i - \mu_0(X_i),\ \ \forall i: T_i = 1 \tag{21}$$

$$\tilde{\tau}_0^i = \mu_1(X_i) - Y_i,\ \ \forall i: T_i = 0 \tag{22}$$

**Step 3: CATE models and weighted fusion.** Two CATE regression models are trained on the treated and control imputations, respectively, and are combined using the propensity score:

$$\tau(x) = e(x) \cdot \hat{\tau}_0(x) + \big(1 - e(x)\big) \cdot \hat{\tau}_1(x) \tag{23}$$

The average treatment effect (ATE) is obtained by averaging the individual treatment effects over the full sample:

$$\tau_{ATE} = \frac{1}{N} \sum_{i=1}^{N} \hat{\tau}(X_i) \tag{24}$$

#### *4.6.4 Sensitivity analysis for unobserved confounding*

To quantify the potential threat of unobserved confounding to the causal conclusions, the E-value

(VanderWeele and Ding, 2017) is used for sensitivity analysis. On the risk-ratio scale, the E-value is defined as the minimum strength of association that an unobserved confounder would need with both the treatment assignment and the outcome in order to fully explain the observed treatment-outcome association. Let the estimated ATE be $\hat{\tau}$ and the control-group baseline risk be $P_0 = E\left[Y \,\middle|\, T = 0\right]$. For $\hat{\tau} > 0$, the risk ratio is approximated by $RR \approx \frac{P_0+\hat{\tau}}{P_0}$. The E-value is computed as:

$$E - value = RR + \sqrt{RR \times (RR - 1)} \tag{25}$$

A larger E-value means a stronger unobserved confounder would be required to overturn the current conclusion, so the conclusion is more robust. An E-value above 1.2 is generally regarded as good resistance to unobserved confounding.

#### *4.6.5 Robustness to the treatment threshold*

The choice of a binarization threshold for a continuous treatment variable is somewhat subjective. To test whether the conclusions depend on a particular threshold, the X-learner procedure is repeated for each treatment variable at four quantile levels (80%, 85%, 90%, and 95%), yielding the range over which the ATE varies with the threshold. If the sign and significance of the ATE remain stable within a reasonable threshold range, the qualitative causal conclusion is considered insensitive to the threshold choice.

# 5. Experiments and analysis

## 5.1 Experimental setup

The experiments were carried out on a mobile workstation running Windows 11, equipped with a 14-core AMD processor, 56 GB of RAM, and an NVIDIA GeForce RTX 4090 Laptop GPU (24 GB of VRAM). The software environment was based on Python 3.10, PyTorch 2.2.0, and CUDA 12.1.

**Table 3.** Search space of the model hyperparameters.

| No. | Parameter | Symbol | Value / range | Role |
|---|---|---|---|---|
| 1 | Learning rate | Lr | 3e-4–1.5e-1 (log) | Step size of weight updates |
| 2 | Batch size | B | 256 | Samples per training step |
| 3 | Number of trees | n_estimators | 200–1200 | Size of the tree ensemble |
| 4 | Number of leaves | num_leaves | 15–125 | Complexity of each LightGBM tree |
| 5 | Tree depth | max_depth | 3–14 | Model depth and overfitting control |
| 6 | Column subsample | colsample | 0.50–1.00 | Feature subsampling to reduce overfitting |
| 7 | Dropout rate | Dropout | 0.1–0.5 | Regularization against overfitting |
| 8 | TCN channels | ch | 64 / 128 / 192 | Multi-scale feature-extraction capacity |

To avoid information leakage, the dataset was split in strict chronological order: the first 70% of the dates

(December 2019 to 25 May 2021, 38,936 samples) were used for training, the middle 15% (26 May to 11 September 2021, about 9,089 samples) for validation, and the last 15% (12 September to 31 December 2021, about 9,212 samples) for testing. The prediction horizon was set to $H \in \{1, 5, 10, 15, 20, 25, 30\}$ days, and the lookback window was fixed at $L = 30$ days after a sensitivity analysis.

The deep models used the AdamW optimizer, a gradient-clipping threshold of 2.0, and an early-stopping patience of six epochs. The enhanced multi-scale TCN had 128 channels, a kernel combination of {5, 7} (retained from the {3, 5, 7} design space after the search), a batch size of 256, and a learning rate of $4.44 \times 10^{-4}$; the GRU had a hidden dimension of 128, a batch size of 256, and a learning rate of $6.19 \times 10^{-4}$. The hyperparameters of the tree models were tuned on the validation set by Bayesian optimization, and the best configuration was determined after 40 iterations; the detailed parameters are listed in Table 3. All classification thresholds were fixed after being chosen on the validation set by maximizing the F1 score, and were then used for test-set evaluation to ensure a fair comparison.

## 5.2 Evaluation metrics

AUC-ROC, the F1 score, accuracy, and recall are used as evaluation metrics.

$$AUC = \int_0^1 TPR \, d(FPR) \tag{26}$$

$$Recall = \frac{TP}{TP+FN} \tag{27}$$

$$F_1 = \frac{2 \times Precision \times Recall}{Precision + Recall} \tag{28}$$

Here TP, FN, and FP are the numbers of true-positive, false-negative, and false-positive samples. All classification thresholds were fixed after being chosen on the validation set by maximizing the F1 score, and were then used for the test set.

Given the special nature of the multi-step fault prediction task, performance is assessed by macro-averaging: the arithmetic mean of each metric over all H horizons, which weighs the performance at every horizon equally.

## 5.3 Ablation study

To clarify the independent contribution of each core component of FGDSE, five comparison variants are designed: TraditionalStacking (no feature partitioning), StackDivNoDeep (the four domain experts only), StackDivGRU (the partitioned experts plus the GRU), StackDivTCN (the partitioned experts plus the enhanced multi-scale TCN), and FGDSE (all modules integrated with horizon-wise dynamic gating). Table 4 and Figure 8 summarize the AUC of the five variants across horizons, together with the decay from H = 1 to H = 30.

The ablation results reveal a clear pattern: the value of feature partitioning is released gradually as the prediction horizon extends. At the short end (H = 1), the five variants cluster tightly within the narrow band of 0.782 to 0.787, a range of only 0.5 percentage point. TraditionalStacking leads at 0.787, StackDivNoDeep follows at 0.782, and FGDSE sits at 0.786, almost indistinguishable from the best. This near-overlapping endpoint matches the objective nature of the data: post faults are strongly first-order autocorrelated at the next-day scale, and the signal-to-noise ratio of this signal is high enough to make architectural choices, such as how features are grouped or whether deep temporal experts are added, negligible. Once the horizon extends beyond ten days, however, the picture changes fundamentally. The autocorrelation of the one-day-lagged signal decays with the horizon, the explanatory contribution of climate, behavioral, and contextual features comes to the fore, and TraditionalStacking, which lacks explicit partitioning, begins to suffer as strong signals suppress weak ones.

**Table 4.** AUC results of the ablation study.

| Model configuration | H=1 | H=5 | H=10 | H=15 | H=20 | H=25 | H=30 | ΔAUC (H=1→30) |
|---|---|---|---|---|---|---|---|---|
| StackDivNoDeep | 0.782 | 0.772 | 0.762 | 0.750 | 0.742 | 0.730 | 0.725 | 0.057 |
| StackDivGRU | 0.783 | 0.774 | 0.766 | 0.754 | 0.744 | 0.731 | 0.727 | 0.056 |
| StackDivTCN | 0.782 | 0.774 | 0.766 | 0.754 | 0.746 | 0.733 | 0.727 | 0.055 |
| TraditionalStacking | **0.787** | 0.777 | 0.770 | 0.758 | 0.749 | 0.735 | 0.730 | 0.057 |
| **FGDSE** | 0.786 | 0.780 | 0.772 | 0.766 | 0.763 | 0.756 | 0.754 | **0.032** |

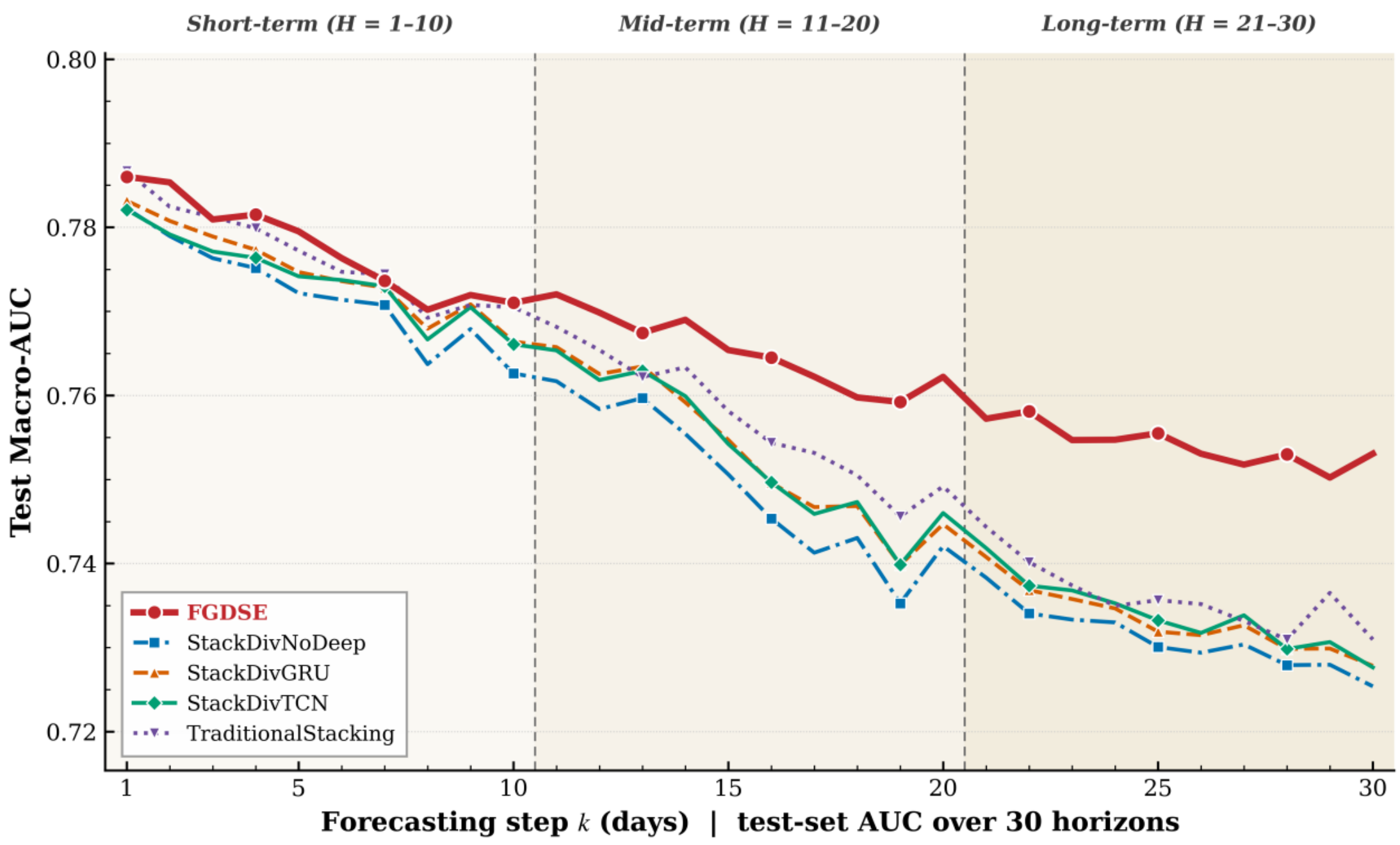

**Figure 8.** Performance comparison of different model configurations.

The quantitative comparison confirms this pattern. From H = 1 to H = 30, the AUC of FGDSE decays by only about 3.2 percentage points (0.786 to 0.754), whereas that of TraditionalStacking decays by 5.7 points (0.787 to 0.730); the accuracy loss of the coarse-concatenation baseline at the monthly horizon is close to 1.8 times that of FGDSE. The single-deep-expert variants StackDivGRU and StackDivTCN (decays of 5.6 and 5.5 points) outperform the pure domain-expert variant but remain clearly below the full FGDSE, which shows that the local temporal memory of the GRU and the multi-granularity trend modeling of the enhanced multi-scale TCN are complementary at long horizons, and only their cooperation can keep suppressing the decay. Notably, at the medium-to-long horizon of H = 15, the AUC of FGDSE exceeds that of TraditionalStacking by about 0.8 point (0.766 versus 0.758), which suggests that the marginal benefit of partitioning is largest in the middle range, in line with the dominant gating-weight transition observed in this interval.

The robustness of recall has more direct operational value. The recall of FGDSE stays near 0.85, decaying by less than 1.2 percentage points from H = 1 to H = 30. This means that even at the long horizon of 30 days, the framework still identifies about 85% of the posts about to fail, providing a reliable basis for monthly spare-part ordering and cross-region resource scheduling.

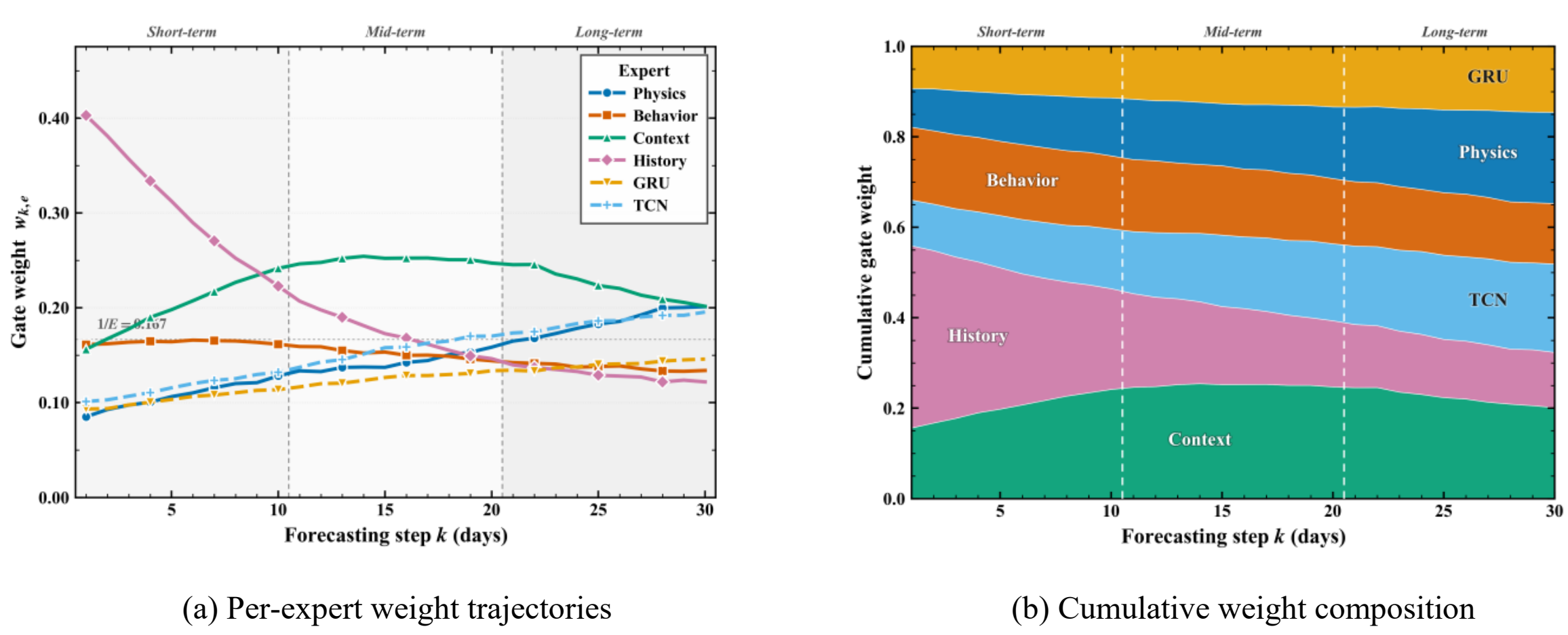

(a) Per-expert weight trajectories

(b) Cumulative weight composition

**Figure 9.** Horizon-adaptive gating behavior of FGDSE.

The evolution of the gating weights along the 30-step horizon (Figure 9) gives a mechanistic explanation for this robustness. In the short-term prediction of one to five days, the historical expert and the behavioral expert jointly hold the largest weight share, taking on the role of intercepting short-term, sudden fault pulses. As the horizon grows, the weights of the physical expert and the enhanced multi-scale TCN expert rise steadily, and the contextual expert also trends slowly upward. This smooth migration of the gating weights matches the physical process by which equipment moves from sudden failure to slow degradation: once sustained heat stress or continuous rain crosses the stage of quantitative accumulation, it gradually turns into real damage such as insulation aging or connector moisture ingress, and the accumulation rate of such damage is clearly no longer driven by the next-day fault history. The model thus hands the lead smoothly from short-term fault experience to long-term cumulative environmental stress, which confirms the physical logic of the charging network moving from sudden failure to slow degradation, and offers the first piece of evidence for the interpretability of the framework.

## 5.4 Comparison with multiple algorithms

### *5.4.1 Baseline algorithms*

To examine the relative advantage of FGDSE systematically, twelve representative algorithms are chosen as baselines. Each baseline is tuned by Bayesian optimization to its best configuration, with detailed parameters given in Appendix A1:

**(1) Persistence.** It uses the fault state of the previous day as the output probability for all future horizons, without any feature learning.

**(2) Logistic regression (Logistic Regression; Hosmer et al., 2013).** With L2 regularization, it serves as a linear baseline to verify the gain of nonlinear modeling over a linear method.

**(3) Random forest (RandomForest; Breiman, 2001).** A Bagging-based ensemble that averages many decision trees, which reduces variance and resists overfitting.

**(4) XGBoost (XGBoost-Full; Chen & Guestrin, 2016).** An end-to-end tree-boosting system trained on the full feature space, without partitioning.

**(5) LightGBM (LightGBM-Full; Ke et al., 2017).** Trained on the full feature space using gradient-based one-side sampling (GOSS) and exclusive feature bundling (EFB).

**(6) LSTM (LSTM-Single; Hochreiter and Schmidhuber, 1997).** A recurrent network that solves the vanishing-gradient problem through gating, and is good at capturing long-term temporal dependence.

**(7) GRU (GRU-Single; Cho et al., 2014).** Comparable in scale to the GRU expert in FGDSE, but run as a standalone model.

**(8) TCN (TCN-Standard; Bai et al., 2018).** A standard single-branch temporal convolutional network with dilated causal convolution and residual connections, without the multi-scale parallel branches, SE channel attention, or attention-based temporal pooling.

**(9) Transformer (Vaswani et al., 2017).** A self-attention network that is good at capturing global long-range dependence in a sequence, included as a baseline for complex temporal modeling.

**(10) Support vector machine (SVM; Cortes and Vapnik, 1995).** A classic method that maps data to a high-dimensional space through a kernel, used here to assess the basic performance of traditional kernel methods on multi-source heterogeneous features.

**(11) CatBoost (Prokhorenkova et al., 2018).** It uses symmetric decision trees and supports categorical variables natively, representing a strong tree ensemble without feature partitioning.

**(12) ARIMA (Box and Jenkins, 1970).** A representative statistical tool for time-series prediction whose forecast relies heavily on the stationarity assumption, used as a baseline for examining the necessity of nonlinear deep modeling.

### *5.4.2 Performance analysis*

After clarifying the internal cooperation of the framework, FGDSE is tested within a comprehensive benchmark of twelve representative baselines. Table 5 reports the multi-step AUC together with the macro-averaged AUC, F1, accuracy, and recall, and Figure 10 visualizes the comparison.

**Table 5.** Performance comparison across algorithms.

| **Model configuration** | **H=1** | **H=5** | **H=10** | **H=15** | **H=20** | **H=25** | **H=30** | **Macro-AUC** | **Macro-F1** | **Macro-Accuracy** | **Macro-Recall** |
|---|---|---|---|---|---|---|---|---|---|---|---|
| **FGDSE** | 0.786 | 0.780 | 0.772 | 0.766 | 0.763 | 0.756 | 0.754 | **0.766** | 0.639 | 0.652 | **0.850** |
| RandomForest | 0.784 | 0.776 | 0.765 | 0.757 | 0.758 | 0.751 | 0.747 | 0.760 | 0.639 | 0.669 | 0.810 |
| LightGBM-Full | 0.788 | 0.781 | 0.765 | 0.758 | 0.751 | 0.740 | 0.739 | 0.759 | 0.639 | 0.671 | 0.806 |
| Transformer | 0.776 | 0.774 | 0.766 | 0.761 | 0.755 | 0.745 | 0.738 | 0.758 | 0.636 | 0.668 | 0.803 |
| XGBoost-Full | 0.786 | 0.776 | 0.762 | 0.748 | 0.744 | 0.736 | 0.742 | 0.754 | 0.633 | 0.665 | 0.798 |
| CatBoost | 0.780 | 0.775 | 0.756 | 0.751 | 0.743 | 0.743 | 0.734 | 0.750 | 0.633 | 0.651 | 0.832 |
| Logistic Regression | 0.786 | 0.768 | 0.751 | 0.722 | 0.720 | 0.712 | 0.729 | 0.737 | 0.622 | 0.645 | 0.813 |
| TCN-Standard | 0.761 | 0.755 | 0.748 | 0.739 | 0.725 | 0.721 | 0.723 | 0.737 | 0.628 | 0.651 | 0.816 |

| GRU-Single | 0.761 | 0.751 | 0.739 | 0.727 | 0.722 | 0.726 | 0.738 | 0.735 | 0.617 | 0.650 | 0.782 |
|---|---|---|---|---|---|---|---|---|---|---|---|
| LSTM-Single | 0.761 | 0.752 | 0.744 | 0.732 | 0.722 | 0.710 | 0.712 | 0.733 | 0.626 | 0.636 | 0.842 |
| ARIMA | 0.764 | 0.751 | 0.743 | 0.728 | 0.717 | 0.702 | 0.692 | 0.727 | 0.613 | 0.625 | 0.821 |
| SVM | 0.714 | 0.700 | 0.670 | 0.650 | 0.642 | 0.630 | 0.642 | 0.660 | 0.534 | 0.621 | 0.602 |
| Persistence | 0.643 | 0.637 | 0.622 | 0.617 | 0.617 | 0.606 | 0.595 | 0.621 | 0.531 | 0.514 | 0.751 |

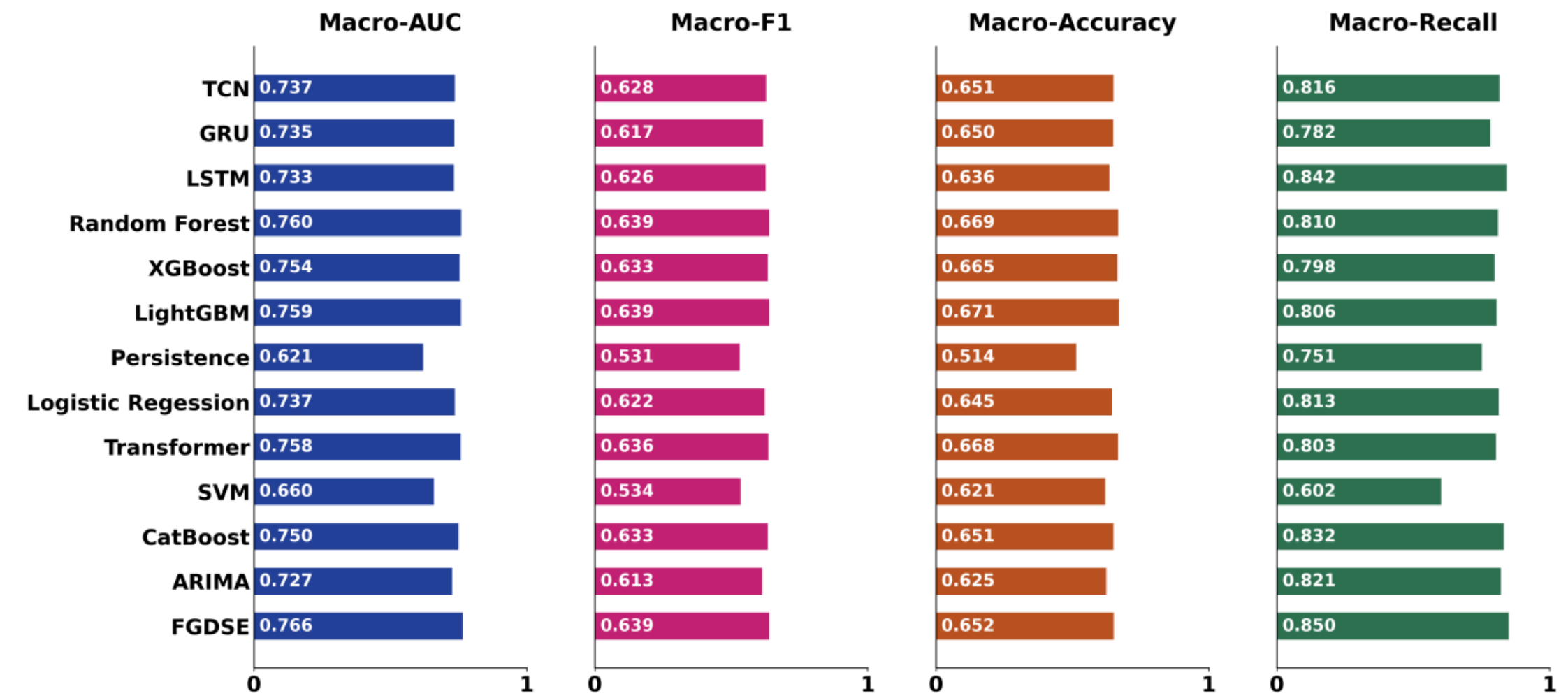

**Figure 10.** Performance comparison across algorithms.

The results can be read along three lines. First, there is clear stratification among model families. Persistence and SVM sit at the bottom (macro-AUC of 0.621 and 0.660), which reflects the limited ability of pure lag transfer and traditional kernel methods on high-dimensional, multi-source heterogeneous features. LSTM, GRU, TCN-Standard, and ARIMA, as general-purpose temporal models, fall in the middle (0.727 to 0.737); their global attention or recurrent structures cannot fully unfold with a limited sample size, and ARIMA, constrained by the stationarity assumption, models the clear nonlinearity and trend drift in the fault sequence poorly. The strong tree models RandomForest and LightGBM-Full lead among the single models (0.760 and 0.759), approaching FGDSE, which confirms the strong nonlinear fitting ability of gradient-boosting trees on tabular heterogeneous features. Transformer ties with the strong tree models at 0.758, which shows that a shallow self-attention network also has comparable fitting potential on monthly operational features.

Second, the advantage of FGDSE is most prominent in the long horizon and in macro-recall. At the very short horizon of H = 1, FGDSE records 0.786, almost overlapping with LightGBM-Full (0.788), Logistic Regression (0.786), and XGBoost (0.786), all within 0.2 percentage point; this reflects that mature single models already have sufficient nonlinear fitting ability when short-term strong autocorrelation dominates, leaving little room for early ensemble gain. As the horizon extends beyond 15 days, however, the gap between FGDSE and the baselines keeps widening: at H = 30, FGDSE reaches 0.754, exceeding RandomForest (0.747), LightGBM-Full (0.739), and CatBoost (0.734) by about 0.7, 1.5, and 2.0 points. The difference in macro-recall carries more operational meaning. FGDSE reaches 0.850, the highest of all models, ahead of the twelve baselines, exceeding RandomForest (0.810) and LightGBM (0.806) by nearly 4.0 to 4.4 points, which on a network of thousands of posts means dozens of additional potential outages intercepted each month. Although LSTM and CatBoost also have high recall (0.842 and 0.832), their macro-AUC is only 0.733 and 0.750, which shows that both tend to issue looser high-risk judgments; their inflated recall comes at the cost of precision, and this differs in essence from FGDSE, which leads in both discrimination and event capture.

Third, stability and interpretability are the dimensions that truly decide engineering value. A single strong tree model essentially crosses multi-source mixed signals in a coarse way; its internal decision path is hard to trace, and it is easily misled by short-term anomalous noise in long-horizon prediction. FGDSE breaks this limit through expert partitioning and dynamic gating. It not only shows better resistance to decay at long horizons but also makes the prediction process structurally traceable: how cumulative physical-environmental stress gradually takes over from short-term fault history can be tracked clearly at the level of the gating weights. For industrial-grade safety warning, a difference in the third decimal of accuracy is often mixed with statistical noise, whereas stable high recall and a transparent decision path are what support real deployment. Based on this engineering reality, this study adopts a clear fault-tolerance stance: it would rather accept a small fluctuation in precision at long horizons than give up the high-recall red line for risk, which is the core value of the partitioned framework in large-scale network risk identification.

## 5.5 SHAP interpretability analysis

To analyze the decision mechanism of FGDSE in next-day prediction, this section uses the dynamic-gating fused probability of the model at H = 1 as a surrogate regression target, trains a LightGBM surrogate, and applies the TreeSHAP algorithm for exact attribution. The coefficient of determination of the surrogate against the original FGDSE output exceeds 0.99, which shows that it reproduces the decision boundary of FGDSE well. The analysis covers all 9,212 test samples and 203 features, and proceeds through four progressive levels: global importance, direction of influence, nonlinear response, and decision path.

FGDSE outputs a 30-day daily risk sequence. Including all 30 horizons in the SHAP analysis would be costly and would also introduce much redundant information. Different horizons correspond to very different operational scenarios. Next-day prediction (H = 1) is directly tied to irreversible, immediate dispatch instructions:

its output triggers next-day work orders and staff assignment, and is the key decision node where predicted probability turns into physical intervention. Weekly and monthly predictions are updated dynamically as data roll forward, and their need for a rigid point-in-time explanation is relatively low. In addition, the analysis of gating-weight evolution in Section 5.3 has already revealed, at the framework level, the macro pattern by which signal dominance migrates as the horizon extends, so this section need not repeat that at the attribution level. In terms of attribution robustness, H = 1 has the highest prediction accuracy among all horizons (AUC = 0.786) and the best fault-autocorrelation signal-to-noise ratio, so its surrogate fits the original output best and is least disturbed by the propagation of approximation error. Weighing the time pressure of operational decisions against the robustness needs of attribution analysis, this section restricts the SHAP analysis to H = 1. To remove the SHAP value sharing caused by collinearity, two features highly correlated with the top feature, FaultMA7d and RecentFault28d, are dropped from the original top 20, and nine mutually independent representative features are kept for the in-depth dependence analysis.

### *5.5.1 Global feature importance structure*

Figure 11 shows the global importance of all features at H = 1. The first fact it reveals is that next-day risk is dominated by fault-history signals, and this dominance is overwhelming. RecentFault14d (mean |SHAP| ≈ 0.055) ranks first, and its contribution is about 1.83 times that of the second-ranked FaultMA14d (≈ 0.030), followed by FaultMA30d (≈ 0.024) and FaultMA7d (≈ 0.018). This ordering reflects the multi-scale structure of post-fault autocorrelation, in which the 14-day window strikes the best balance between sensitivity to recent dynamics and statistical stability, and is the best single estimator of recent fault density. Note that RecentFault14d and FaultMA14d both correspond to the 14-day window, but the former counts cumulative faults and reflects evidence of a state switch, while the latter computes a rolling mean fault intensity and estimates the baseline fault rate. Both enter the top two at H = 1 with similar magnitudes, which shows that the model uses two complementary logics in next-day prediction: a state-switch signal and a baseline-rate estimate. The fault-history features as a whole, including rolling means, cumulative counts, baseline rate, and the per-user ratio, account for about 90% of the total |SHAP| of the top 20, establishing their overwhelming dominance in next-day risk attribution.

The cross-family interaction features MWI×Fault_t-1 (≈ 0.007) and Temp×Fault_t-1 (≈ 0.006) enter the top 10 right after, with magnitudes comparable to first-order aggregate statistics such as FaultRate (≈ 0.006), together accounting for about 9%. This result directly supports the basic premise of FGDSE—feature partitioning followed by fusion. There is a non-additive coupling between meteorological stress and fault history across feature families, and this coupling cannot be captured by any single expert alone; it must be modeled explicitly through cross-family interaction terms. In other words, the SHAP attribution provides independent evidence, at the feature level, for the necessity of the multi-expert fusion architecture.

Compared with the fault-history family, the mean |SHAP| of pure climate features and user-behavior

features does not exceed 0.005, about one order of magnitude below the top feature, and Location ranks only 16th at 0.002. This does not mean that climate, behavioral, and spatial signals are worthless; rather, at the very short horizon the signal-to-noise ratio of the first-order fault-autocorrelation signal is far higher than that of other signal classes, and the action of slow-varying factors needs a longer accumulation time before it shows up in the prediction error. This phenomenon mutually corroborates the gating-weight evolution reported in Section 5.3: at H = 1 the gate leans clearly toward the historical and behavioral experts, while the weight of the physical expert rises only after H ≥ 10. The weight migration at the architecture level and the attribution strength at the feature level reach a consistent conclusion on the same data, forming an internal cross-validation of the interpretability of FGDSE.

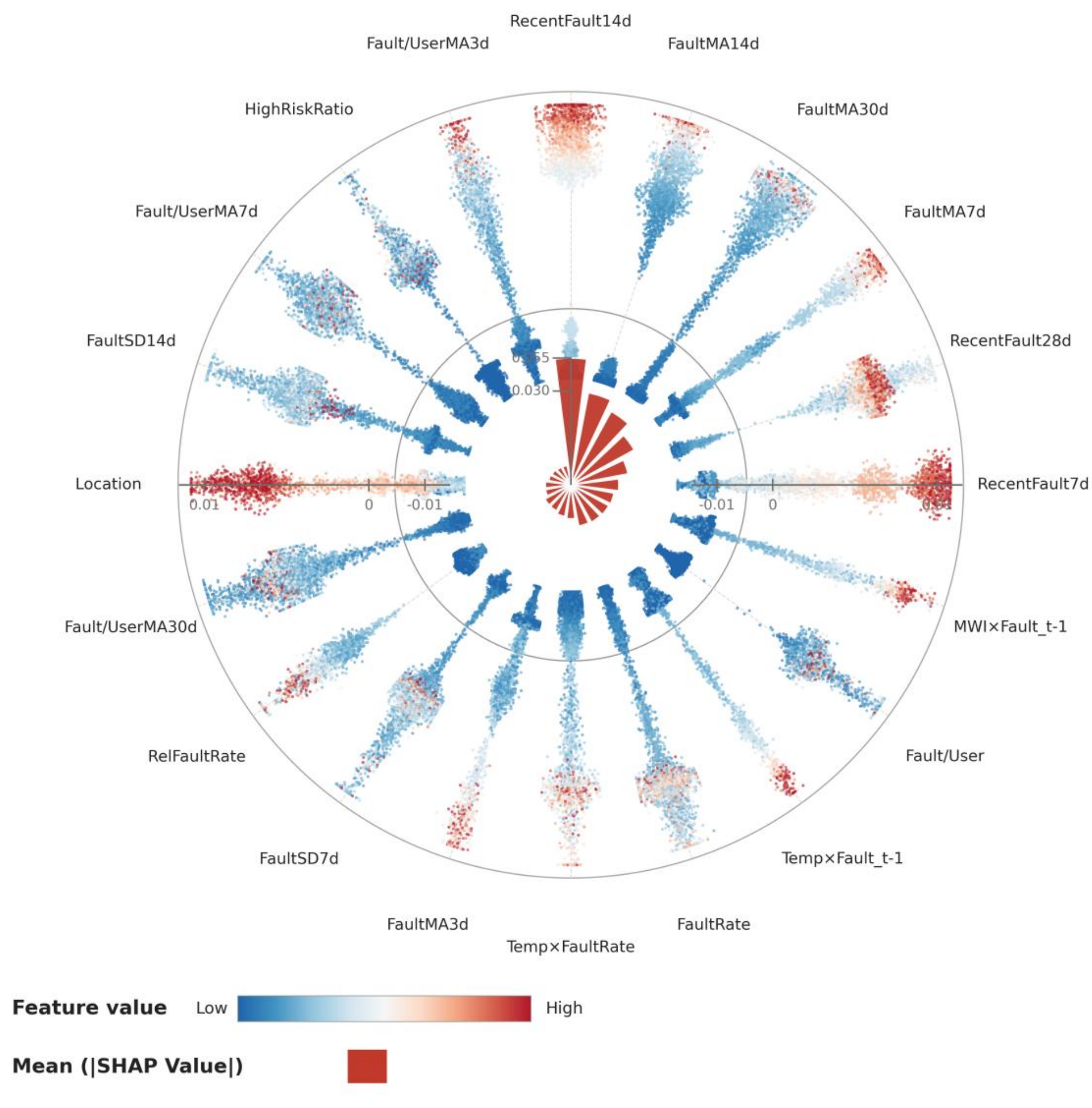

**Figure 11.** SHAP polar beeswarm plot of global feature importance.

Direction and heterogeneity can be read jointly from the color distribution in the outer ring of Figure 11 and the beeswarm plot in Figure 12. The high-value samples of FaultMA30d (red) concentrate in the positive

SHAP region, with some contributions reaching about +0.21; the low-value samples (blue) extend into the negative region down to about −0.15, an almost symmetric distribution. This shows that posts with dense recent faults are clearly raised in next-day risk, while posts long free of faults are actively lowered, which is exactly the core computational advantage of preventive maintenance over after-the-fact response: it keeps limited labor and spare parts from being spread over targets that look calm but are not truly dangerous. The response of FaultRate is clearly heterogeneous: at the same historical fault rate, the SHAP value spreads from +0.02 to near zero or even slightly negative, which shows that the same baseline rate does not contribute equally to next-day risk in different operational contexts. If an operations platform graded posts only into three coarse levels by historical fault rate, it would systematically underestimate highly sensitive posts in special states, such as aging equipment that has served more than 24 months, that bears high peak load, or that has just gone through a rainy season. This heterogeneity is not statistical noise but a real, identifiable difference among subgroups.

### *5.5.2 Beeswarm plot: joint reading of feature-value direction and SHAP effect*

The beeswarm plot (Figure 12) visualizes the correspondence between feature-value direction and model response. Take RecentFault14d, the most important feature at H = 1, as an example. Its beeswarm distribution shows a clear symmetric structure: the red points (high feature values, that is, posts with high recent cumulative faults) gather at the right end of the positive SHAP region, with some samples reaching +0.09, while the blue points (low feature values, that is, posts with no recent faults) extend leftward to about −0.06. This two-sided distribution carries a clear physical meaning: posts with dense recent faults have their risk markedly raised by the model, while posts long free of faults have their risk actively lowered, which is the computational embodiment of the idea that preventive maintenance is better than after-the-fact response.

The beeswarm distributions of the three rolling means—FaultMA14d, FaultMA30d, and FaultMA7d—share a common asymmetric shape. The low-value samples gather into a dense cluster in the negative SHAP region, roughly −0.025 to −0.040, with a narrow span; the high-value samples spread into the positive region, reaching +0.04 to +0.05, with a wider span. This long positive tail and short negative tail show that the safe state represented by a low fault history is relatively uniform, whereas a high fault history may correspond to several risk levels, so the model has more room for differentiation when adjusting upward.

The beeswarm distribution of FaultRate shows a similar but more dispersed pattern: the SHAP values of posts with a high static fault rate concentrate in the positive region but spread more widely. This shows that the same fault-rate level does not contribute equally to the risk estimate under different operational contexts, such as equipment age, region, or load intensity, and that rating by historical fault rate alone may underestimate the true risk of highly sensitive posts in particular operational states.

The cross-family interaction features MWI×Fault_t-1 and Temp×Fault_t-1 show a typical threshold-activation shape. A dense cluster in the negative SHAP region corresponds to the stable negative contribution when the interaction term is low, while the positive region extends to +0.02 to +0.03, corresponding to positive activation when it is high. This shape physically corresponds to the synergy by which meteorological stress amplifies the baseline fault rate: high humidity alone or a prior-day fault alone is not enough to trigger a clear upward adjustment, and only when both are high is the effect activated. This directly verifies, at the level of feature effects, the necessity of introducing cross-family interaction terms in FGDSE.

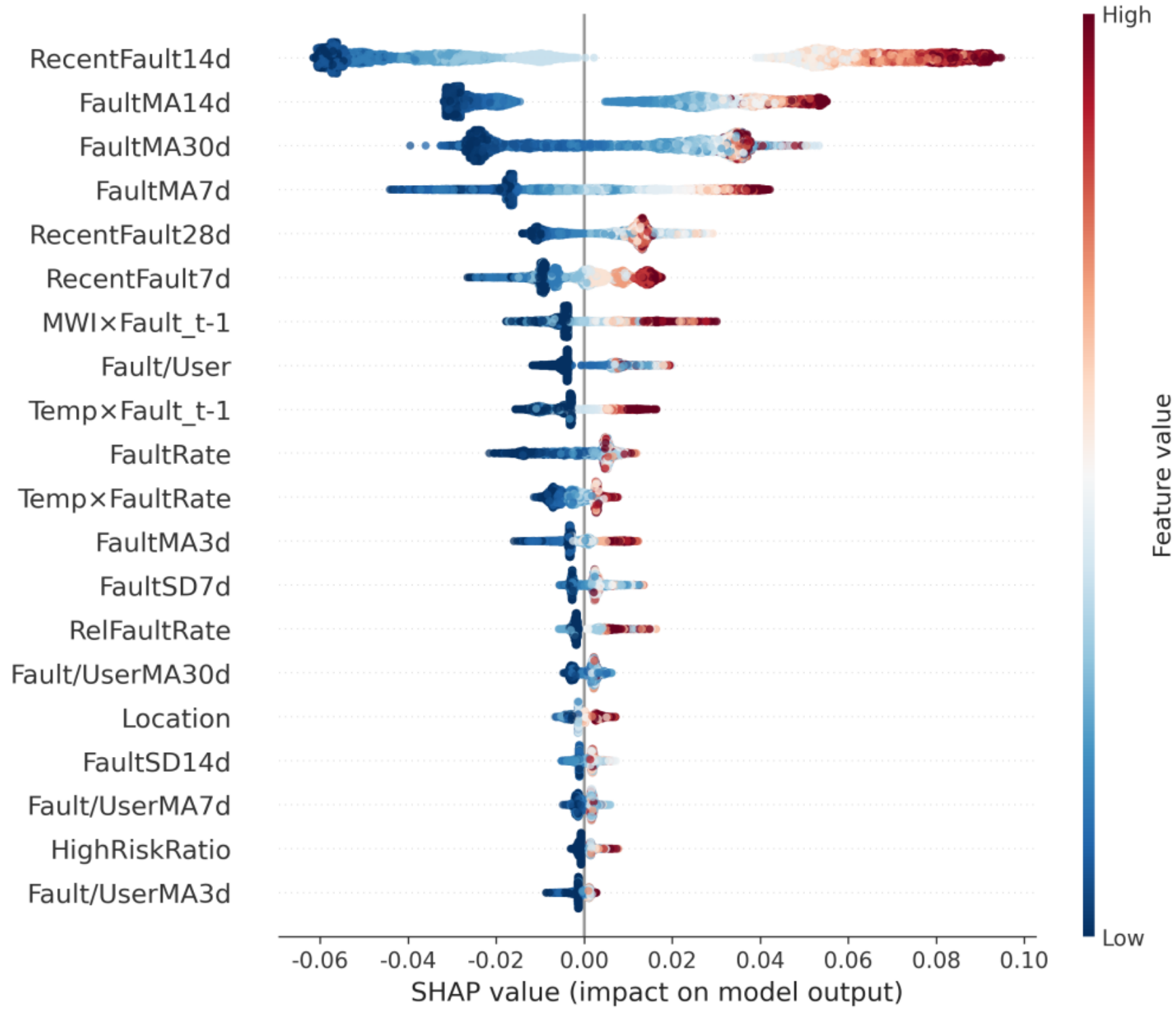

**Figure 12.** SHAP beeswarm plot.

### *5.5.3 Dependence plots: nonlinear response curves between features and risk*

The dependence plots (Figure 13a–i) characterize how risk responds as a feature value changes continuously, and their primary feature is a threshold-saturation response. Take FaultMA14d as an example: as the 14-day rolling mean rises from zero to about one or two, the SHAP value climbs quickly from about −0.025 to about +0.04; thereafter, even as the mean keeps rising to six or even 14, the SHAP increment hardly changes. FaultMA30d shows the same shape, with a similar threshold position but a slightly lower saturation ceiling, which reflects that the long window dulls extreme values more strongly than the medium window. This saturation structure shows that the model judges fault-accumulation intensity with a clear state-transition property: once a post crosses a certain accumulation-density boundary, it is judged to be in a high-risk state, and further accumulation has a diminishing marginal effect on priority. This structure maps directly to an executable maintenance rule: any post that reaches an average of more than one fault per day within the 14-day window enters the next-day priority list, and this threshold rule covers the vast majority of objects the model identifies as high-risk.

The dependence plots of RecentFault14d and RecentFault7d show two different shapes: a monotonic continuous curve and a discrete staircase. The response curve of RecentFault14d rises almost linearly: as the feature value increases from zero to 1.0, the SHAP climbs smoothly from −0.06 to +0.09, with no clear turning point or saturation plateau, so the model treats it as a continuous risk-dose indicator effective across the whole dynamic range, which indirectly explains why it receives the highest importance. Because RecentFault7d is essentially a finite count, its horizontal axis shows seven or eight discrete steps, each corresponding to the marginal effect of one more fault, with a step height of about 0.005 to 0.010. The dependence plot of Location likewise shows a typical categorical-binning structure, with the central SHAP value under different station codes distributed between −0.005 and +0.006, which reflects the inherent role of spatial heterogeneity as a baseline risk offset; this structure suggests that Location acts both as an independent baseline risk indicator and as a modulator of the fault-history effect.

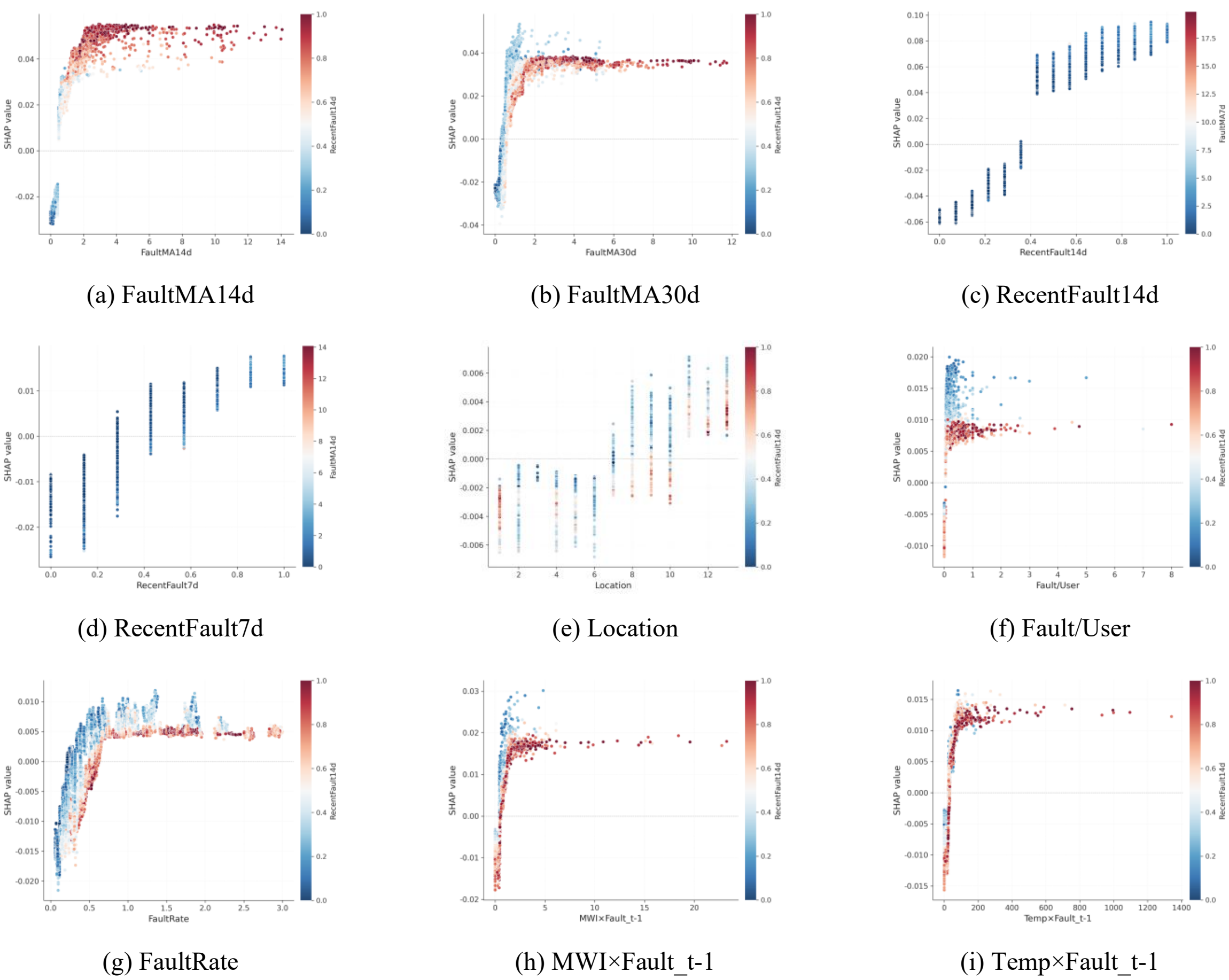

(a) FaultMA14d (b) FaultMA30d (c) RecentFault14d

(d) RecentFault7d (e) Location (f) Fault/User

(g) FaultRate (h) MWI×Fault_t-1 (i) Temp×Fault_t-1

**Figure 13.** SHAP dependence plots.

The dependence plots of Fault/User and FaultRate show two special shapes related to the baseline rate. At the same value of Fault/User, the SHAP values fall into two clusters, one positive at about +0.008 and one negative at about −0.012. In the mathematical sense of SHAP, this vertical spread is not a bimodality of the feature itself but comes from the modulation of an interaction effect. The samples with a recent high fault frequency gather in the positive SHAP region, while those free of recent faults sink into the negative region; in other words, the model does not treat the per-user fault rate as independent evidence but judges it after coupling with the recent fault state. The dependence plot of FaultRate is highly similar in shape to that of Fault/User. The median SHAP trajectory rises slowly with the feature value, but at the same FaultRate level the scatter spans about 0.02, and the colors form a clear layering with red above and blue below. This shows that the model does not view the historical fault rate independently and linearly but reads it in strong coupling with the recent fault state: posts with a high historical rate that are still faulting recently are raised markedly, while posts with a high historical rate that have stabilized recently are lowered.

The dependence plot of MWI×Fault_t-1 shows an almost flat response plateau when the product is below one, with the SHAP staying near zero; once the product exceeds one, the SHAP starts to rise, entering a saturation region at about five, with a maximum contribution of about +0.02. Temp×Fault_t-1 shows a highly similar activation threshold and saturation ceiling, which further shows that the coupling between meteorological stress and fault history is not an isolated case but a stable pattern across several specific meteorological variables. This slow-start shape shows that there is an activation threshold between meteorological stress and the prior-day fault count: high temperature alone or a prior-day fault alone is not enough to trigger a clear risk increase, and the coupling effect is activated only when both are high. The corresponding executable maintenance strategy is to give priority, during a spell of sustained heat or high humidity, to posts that reported a fault the day before, rather than spreading inspection resources evenly over all equipment in use.

#### *5.5.4 Decision plot: class comparison of individual prediction paths*

The decision plot (Figure 14) verifies the discrimination quality of the model at the individual level. At H = 1, the path behavior of 30 randomly drawn fault samples differs sharply from that of 30 normal samples: the fault samples start from a baseline of about 0.42 and, under the cumulative action of the fault-history features, gradually shift right to the 0.55–0.75 range, while the normal samples shift the other way to 0.18–0.35, forming a clear scissor-like separation.

The main contribution to the separation concentrates in the first six to eight features, led by RecentFault14d, FaultMA14d, FaultMA30d, and FaultMA7d, after which the remaining features only adjust the path slightly, in line with the concentrated importance shown by the global attribution.

This separation means that every post flagged by FGDSE as high-risk can be traced to the specific feature names and amounts that contribute most, forming a dual evidence chain in which the risk score and the attribution explanation are placed side by side. In a real dispatch workflow, this lets every scheduling instruction carry a short risk summary—for example, noting that the post has a high 14-day cumulative fault count, that its 30-day mean has exceeded the warning threshold, and that the day falls in a spell of sustained high humidity, so the next-day fault probability is judged high and priority handling is advised for the shift on duty. A traceable basis for decisions is the key difference between FGDSE and a black-box deep model, and a precondition for its adoption as an operations tool.

Taken together, the interpretability of FGDSE can be summarized at three mutually corroborating levels. At the feature level, the short term is overwhelmingly dominated by fault-history features, followed by cross-family interaction terms, while the long term is taken over by climate trends and contextual factors. At the individual level, high-risk and normal samples separate systematically along their cumulative paths, with the first six to eight features completing the main discrimination. At the framework level, the fault-history dominance shown by SHAP at H = 1 is highly consistent with the short-term gating-weight evolution in Section 5.3, and the threshold-activation and synergistic-amplification mechanisms identified in Section 5.5.3 echo the stable advantage of FGDSE over single-expert baselines in Section 5.4. The evidence at the three levels together forms an internally consistent chain for the interpretability of FGDSE, and provides a quantitative, mechanism-based anchor for later operational decisions.

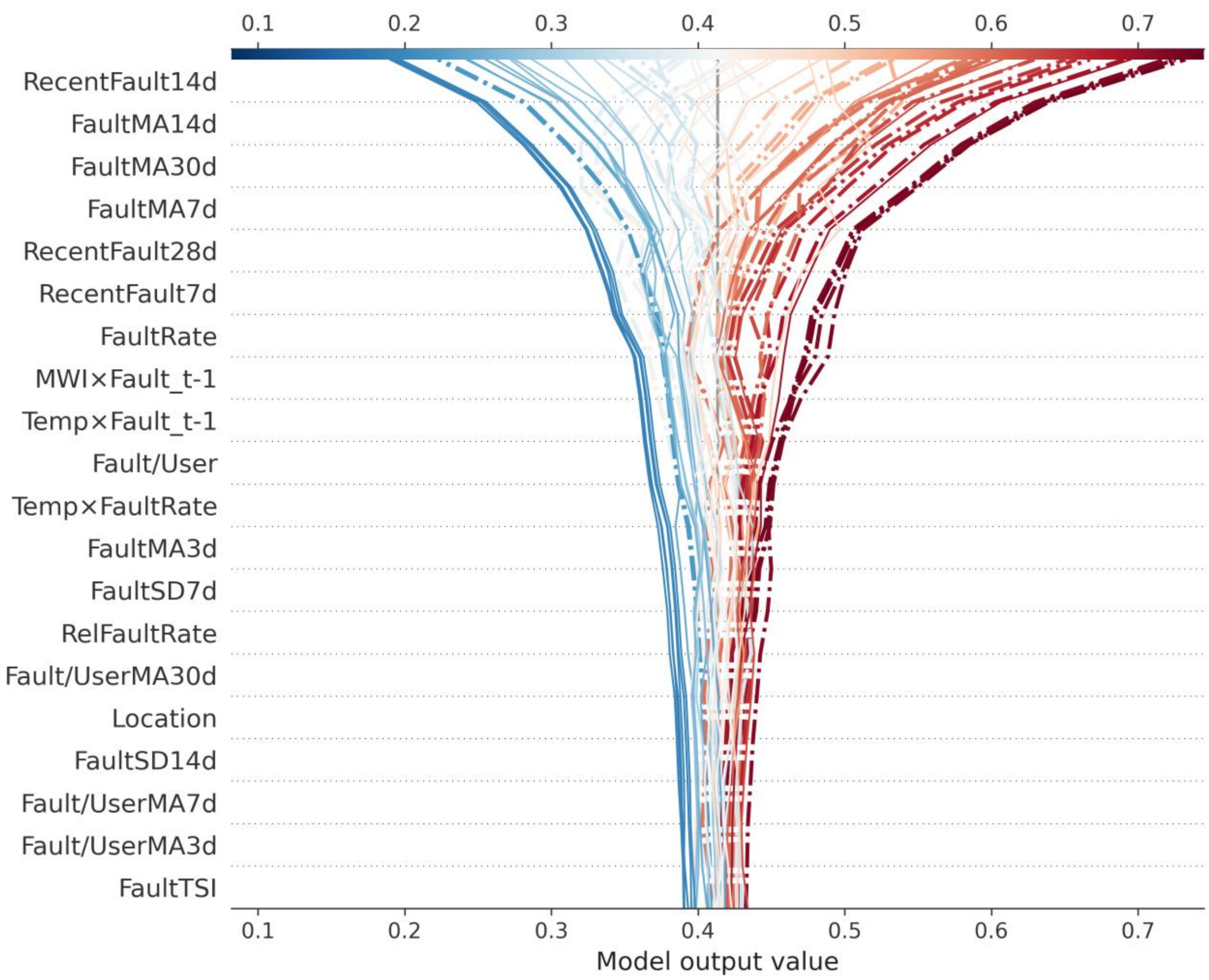

**Figure 14.** SHAP decision plot.

## 5.6 Causal inference analysis with the X-learner

To move beyond the correlational limit of SHAP attribution, this section introduces the X-learner causal inference framework. Based on 57,237 observational records, it estimates the conditional average treatment effect (CATE) of four core treatment variables across horizons. The treatment threshold is set adaptively at the 85th or 90th percentile of the meteorological and operational indicators, and the full statistical estimates are given in Table 6.

### *5.6.1 Diagnostics and robustness checks for the causal assumptions*

Before interpreting specific effects, the basic premises on which the causal inference relies must be checked. Figure 15(a) gives the probability density of the propensity score under each treatment variable. The distribution curves of the treatment and control groups overlap closely in the core probability range, and the overlap ratio exceeds 0.95 in all experiments. This shows that the observational data have an adequate common support in the multidimensional covariate space, so the treatment and control groups are well comparable. The covariate-balance diagnostics confirm this further. After propensity-score weighting, the standardized mean difference (SMD) of most experiments falls below 0.02, which shows that the covariate distributions have become balanced.

The case of extreme-heat treatment, however, needs separate comment. After weighting, its maximum SMD reaches 0.308, clearly higher than that of the other experiments. The deviation arises because the covariate set naturally lacks a feature that is both correlated with temperature and not highly collinear with the treatment variable. Including a temperature-derived feature would make the propensity-score estimation unstable through severe multicollinearity. The E-value is further used to check the robustness of the extreme-heat effect to unobserved confounding, giving an E-value of 1.23. This means that any unobserved confounder seeking to fully explain the observed causal association would need to raise both the treatment-assignment probability and the fault risk by at least a factor of 1.23. The extreme-heat effect therefore remains fairly robust to unobserved confounding.

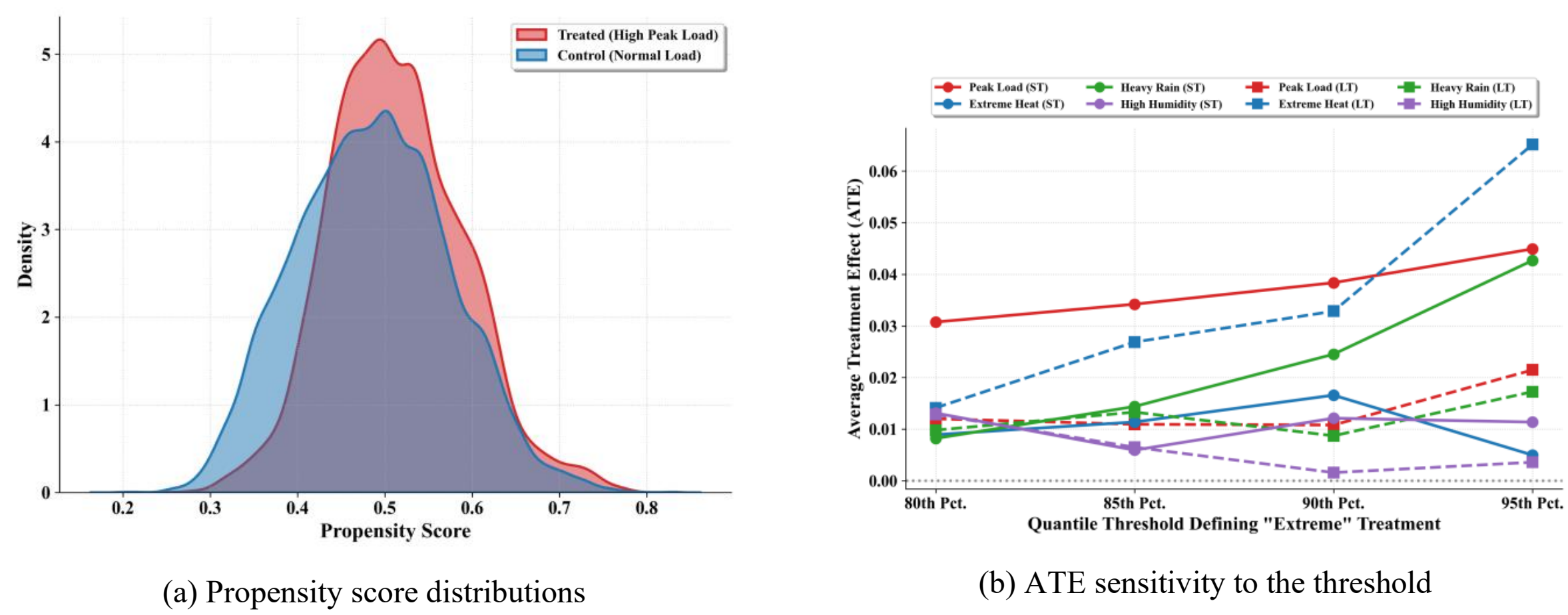

(a) Propensity score distributions

(b) ATE sensitivity to the threshold

**Figure 15.** Robustness checks for the causal inference.

In addition, the threshold-sensitivity analysis supports the reliability of the core conclusions. Figure 15(b) shows how the ATE of each treatment variable changes as the quantile threshold is tightened from 80% to 95%. The positive effect of high peak load on short-term faults stays within 0.0308 to 0.0449, and the positive effect of extreme heat on long-term faults stays within 0.0141 to 0.0652. The sign of neither effect ever reverses as the threshold changes. Notably, the long-horizon ATE of extreme heat jumps from 0.0329 to 0.0652—nearly doubling—over the narrow interval where the threshold moves from 90% to 95%. This nonlinear jump suggests a possible threshold effect in heat-induced damage; however, since the jump appears only in the narrow 90%–95% interval and the sample size falls accordingly, this reading needs confirmation on a larger sample. By comparison, the short-term effect of heavy rain has a relatively wide range, from 0.0082 to 0.0427, though its sign stays positive throughout.

### *5.6.2 Heterogeneous causal effects and their mechanisms*

Among the four core variables, high peak load and heavy rain show an acute short-term positive effect (short-horizon ATE of about +0.0342 and +0.0245), which shows that, after controlling for all other covariates, these two physical shocks raise the fault probability of posts significantly in the short term.

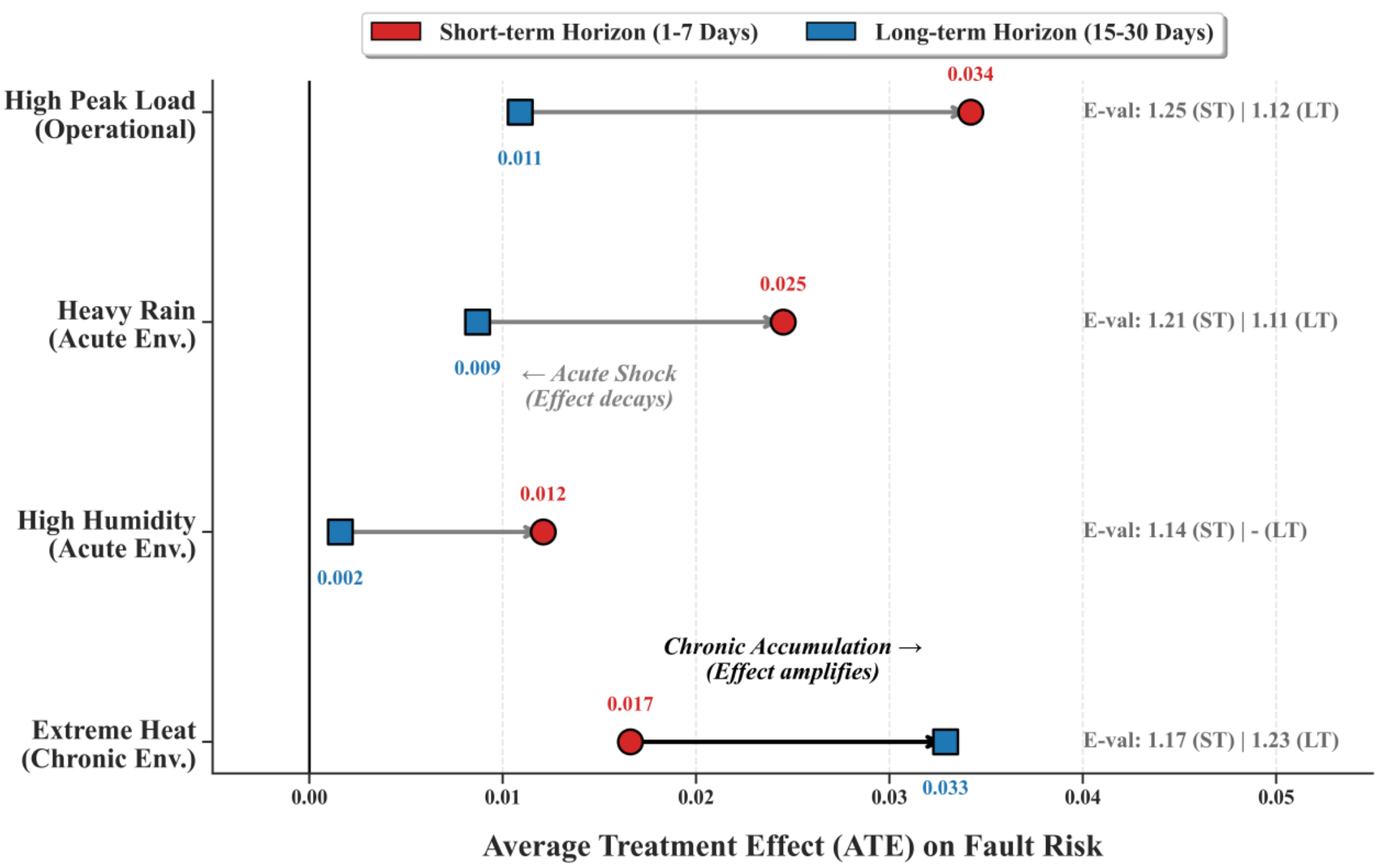

**Figure 16.** Spatiotemporal evolution of the X-learner causal effects.

This conclusion has direct mechanistic support: high-intensity load easily causes instantaneous thermal sintering of contactors, while rainwater ingress easily short-circuits internal components, both of which are direct triggers that breach the safety threshold of the equipment. In terms of heterogeneity, the vulnerable subgroup whose short-term fault probability rises by more than 5 percentage points (CATE > 5pp) under high peak load reaches 31.3%, forming the high-risk group that should receive priority peak-shaving and load limiting. As the horizon extends to 15 to 30 days, however, the causal effects of both decay sharply (long-horizon ATE down to about +0.0109 and +0.0087), and the high-risk subgroup share falls steeply to 12.0% and 13.8%. This decay not only defines the action boundary of acute stress but also shows, from the data, that instantaneous overload or short-term water ingress lacks a long-term cross-horizon memory effect.

In sharp contrast to the acute shocks, extreme heat is the only variable that shows a cross-horizon amplification effect (short-horizon ATE of about +0.0166, rising to about +0.0329 at the long horizon). Figure 16 shows how the ATE of the four treatment variables evolves over the 1-to-30-day horizon, as a reference for comparing the acute-shock and cross-horizon-amplification mechanisms. This asynchronous phenomenon also has solid mechanistic support: the thermal-fatigue cycling rate of the power-conversion module accelerates exponentially with junction temperature, the contact resistance at the connector interface becomes unstable due to metal creep, and the aging rate of insulation material rises markedly with temperature following the Arrhenius relationship. The CATE distribution further reveals the heterogeneity of heat damage: at the long horizon, as many as 30.4% of posts have a fault-probability rise above 5 percentage points, forming the absolute priority for summer preventive maintenance; another group of posts has a CATE that is not significant or even negative, which may be explained by operators having already applied active inspection or power derating to them during hot periods, an unobserved supply-side response that appears in the data as a reverse drop in the fault rate of some posts.

**Table 6.** Estimation results of the X-learner.

| Variable | Horizon | Threshold | n_treated | ATE | 95% CI | CATE>5pp | E-value | ATE range (80%→95%) |
|---|---|---|---|---|---|---|---|---|
| High Peak Load | 1–7 days | >85% | 8586 | 0.0342 | [0.034, 0.035] | 31.3% | 1.25 | [0.0308, 0.0449] |
| High Peak Load | 15–30 days | >85% | 8586 | 0.0109 | [0.010, 0.011] | 12.0% | 1.12 | [0.0108, 0.0215] |
| Extreme Heat | 1–7 days | >90% | 5709 | 0.0166 | [0.016, 0.018] | 37.2% | 1.17 | [0.0050, 0.0166] |
| Extreme Heat | 15–30 days | >90% | 5709 | 0.0329 | [0.032, 0.034] | 30.4% | 1.23 | [0.0141, 0.0652] |
| Heavy Rain | 1–7 days | >90% | 5711 | 0.0245 | [0.024, 0.025] | 29.9% | 1.21 | [0.0082, 0.0427] |
| Heavy Rain | 15–30 days | >90% | 5711 | 0.0087 | [0.008, 0.009] | 13.8% | 1.11 | [0.0087, 0.0172] |
| High Humidity | 1–7 days | >90% | 5474 | 0.0121 | [0.011, 0.013] | 24.3% | 1.14 | [0.0059, 0.0131] |
| High Humidity | 15–30 days | >90% | 5474 | 0.0016 | [0.001, 0.002] | 11.4% | 1.04 | [0.0016, 0.0128] |

The high-humidity variable deserves separate discussion: its long-horizon ATE is close to zero (about +0.0016). This counterintuitive result—a harsh environment showing no severe damage—can be understood through a dual mechanism. The first is demand-side behavioral adaptation: bad weather suppresses users' willingness to travel, the daily number of transactions per post drops, and the mechanical plugging and power cycling of connectors fall accordingly. The second is supply-side maintenance response: defensive inspection launched by operators during weather warnings actively suppresses the fault rate. However, the weak average effect at the long horizon hides a key heterogeneous difference: the share of posts with a CATE above 5 percentage points still reaches 24.3% at the short horizon and 11.4% at the long horizon, which shows that the network still contains specific subgroups extremely sensitive to humidity, whose high risk comes mainly from slow moisture ingress caused by seal aging. The analysis above points to a core proposition: the influence of climate and operational factors on charging-facility fault risk must be built on post-level CATE estimation, so that maintenance resources can be allocated precisely and differentially.

## 5.7 Research implications

Drawing together the experimental analysis above, this study offers important implications in five respects.

**(1) Dual-window forecasting supports the practical adoption of preventive maintenance, which in turn raises charging-network availability, eases range anxiety, and protects the emission-reduction benefits of replacing fuel vehicles with EVs.** At the 30-day horizon, FGDSE still keeps a macro-recall of 0.85 with a decay below 1.2 percentage points, so operations teams can lock in high-risk posts a month ahead. Specifically, operators are advised to output short-term urgent work orders and a monthly maintenance plan over two windows of 7 and 30 days: the higher-recall warnings within 7 days directly trigger inspection dispatch, while the 30-day forecast is used to allocate cross-region spare-part stock and labor. This dual-window scheduling shifts the management logic from passive response to active intervention and lowers the daily revenue loss per station caused by fault downtime. From the perspective of transport-environment governance, higher availability not only adds usable post-sessions each day and eases charging congestion in dense cities, but also reduces detours, queuing, and emergency-repair trips spent looking for an available post, which indirectly lowers transport energy use and emissions. More importantly, reliable charging service lowers users' perceived uncertainty about EV travel and sustains the long-term behavioral basis for replacing fuel vehicles with EVs.

**(2) Climate-adaptation retrofits should be invested by post-level CATE grading, avoiding a one-size-fits-all strategy.** The long-horizon X-learner analysis reveals that about 30.4% of posts have a significant rise in fault probability under extreme heat (CATE > 5 percentage points), and the short-horizon share is even higher at 37.2%. Operators are therefore advised to build a climate-vulnerability scorecard for posts: for stations with heat-sensitive posts whose CATE exceeds 5 percentage points, give priority to adding active cooling ducts or sun-shading structures and schedule a dedicated thermal inspection before each summer; for posts with a negative CATE, distill their operational experience, such as derated operation and timed ventilation, into a reproducible standard operating procedure and promote it to highly vulnerable posts. For the humidity subgroup with a CATE above 5 percentage points (24.3% at the short horizon and 11.4% at the long horizon), connector seal replacement and cabinet anti-moisture coating checks should be completed before the rainy season. This

grading strategy directs a limited maintenance budget precisely to where the marginal benefit is highest.

**(3) The horizon-dependent shift in signal dominance provides a theoretical basis for a step-by-step maintenance strategy.** SHAP and the gating weights jointly show that short-term prediction is dominated by fault history, while long-term prediction is taken over by climate and context. Accordingly, a short-term warning system should use recent 14-day fault-density indicators as its core trigger, issuing a red alert once the daily rolling mean of faults exceeds one per day; quarterly or annual maintenance plans, by contrast, should focus more on the long-term trends of equipment service life, the regional climate-exposure index, and load concentration during peak hours, and set the annual priority for preventive inspection on that basis. Decoupling the two indicator systems prevents short-term noise from disturbing long-term decisions.

**(4) The equipment-health information embedded in charging behavior offers a new angle for refined operation.** The SHAP analysis shows that user repeat-use rate and peak-hour concentration are significantly associated with fault risk, corresponding to the cumulative effects of connector mechanical wear and thermal fatigue. Operators are advised to use time-of-use pricing to shift charging load toward off-peak hours—for example, by deepening off-peak price discounts, or by introducing a queuing-and-reservation mechanism during peak hours at stations with high-fault posts—to flatten the instantaneous load peak and slow the accumulation of equipment fatigue. The share of zero-power charging can also be added to the equipment-health dashboard: when this indicator exceeds 1.5 times its historical mean for three consecutive days, a connector-contact inspection order is triggered automatically, achieving early and active discovery of hidden faults. From a demand-response perspective, the load redistribution guided by time-of-use pricing not only slows equipment-fatigue accumulation but also reduces the grid's reserve-capacity need through peak-shaving, indirectly lowering marginal generation start-ups and producing an implicit carbon-reduction benefit.

**(5) Causal inference provides a quantitative basis for ex-ante policy evaluation and low-carbon transport governance.** Unlike traditional correlation analysis, the CATE estimated by the X-learner can answer the counterfactual question of how much the fault probability changes after a specific intervention is applied to a specific post. This ability can be extended to evaluate the expected effect of concrete policies such as derated operation, periodic seal inspection, and intensified summer inspection. For example, based on the post-level heat CATE, an operator can estimate the fault-rate reduction after a 10% power derating for highly vulnerable posts in summer, and carry out a cost-benefit analysis against the charging-energy loss caused by derating, providing a scientific basis for climate-adaptation investment in charging infrastructure. From a broader low-carbon-mobility view, higher charging-network availability directly eases users' range anxiety and advances the replacement of fuel vehicles with EVs, so the environmental benefit of preventive maintenance lies not only in reducing emergency-repair trips but also in sustaining a stable supply of low-carbon travel services.

## 6. Conclusion

This paper addresses multi-step fault risk prediction for electric vehicle charging facilities—a scientific problem of considerable practical value and theoretical significance—and proposes the FGDSE framework, which jointly advances three dimensions: prediction accuracy, multi-step sequence capability, and interpretable causal inference.

(1) The value of the feature-partitioning strategy is released markedly as the prediction horizon lengthens. The AUC of the FGDSE dynamic gating framework exceeds that of all baselines at $H \geq 10$; at very short horizons it is essentially on par with the best single model, and its advantage over TraditionalStacking widens steadily as the horizon extends. From $H = 1$ to $H = 30$, the AUC of FGDSE decays by only 3.2 percentage points (0.786 to 0.754), well below the 5.7 percentage points of the traditional stacked ensemble (0.787 to 0.730). This pattern confirms the theoretical expectation that feature partitioning has limited effect in short-term prediction but is indispensable in long-term prediction, and it reveals a transition in signal dominance: fault-historical signals dominate short-term prediction, while climate and contextual signals gradually take over at long horizons.

(2) The recall of the framework remains highly robust at very long horizons. Even at the 30-day horizon, the framework still identifies about 85% of the charging posts about to fail, with a decay of less than 1.2 percentage points. This property directly supports real operational scenarios such as monthly resource pre-allocation and spare-part pre-ordering.

(3) Causal inference extends the prediction system toward actionable intervention decisions. The X-learner analysis identifies extreme heat as the only variable whose positive causal effect amplifies as the prediction horizon extends (from a short-horizon ATE of about +0.017 to a long-horizon ATE of about +0.033); at short horizons about 37.2% of charging posts show a fault-probability increase of more than 5 percentage points during extreme heat, and at long horizons this share is 30.4%. The subgroup sensitive to high humidity accounts for 24.3% at short horizons and 11.4% at long horizons. The heterogeneity of post-level CATE provides a quantifiable basis for allocating differentiated maintenance investment and extends the prediction system from a risk score into executable causal decision support.

This study has several limitations that warrant deeper investigation. First, the predictive capability of the current framework relies on transaction logs and meteorological data; introducing high-frequency physical quantities, such as the internal temperature of charging posts and connector plug-in torque, could further improve short-term prediction and would provide a more direct data basis for embedding physical prior information. Second, the causal estimation of the X-learner relies on the ignorability assumption, and operational response behavior—an unobserved variable—poses a substantial threat to this assumption; quasi-experimental evaluation based on natural experiments is an important direction for strengthening the credibility of the causal findings. Third, the current framework treats each horizon as an independent classification target and does not explicitly model the conditional dependence between adjacent horizons; generative models for sequence labels may further improve accuracy by capturing cross-horizon correlation. In addition, incorporating the geographic proximity between stations into the ensemble through a graph neural network, so as to capture spatial-clustering effects, is a natural extension of FGDSE from the single-post level to network-level risk management.

# Appendix

Table A1. Optimal hyperparameter settings for baseline models.

| Model | Hyperparameter configuration | Notes |
|---|---|---|
| **Persistence** | strategy = yesterday_fault_binary | Naive baseline that propagates the previous day's fault state to all 30 horizons; no learnable parameter. |
| **Logistic Regression** | C = 0.374, penalty = l2, solver = liblinear, class_weight = balanced, max_iter = 1000, random_state = 42 | L2-regularized linear classifier; class_weight = balanced addresses the fault / non-fault imbalance. |
| **RandomForest** | n_estimators = 150, max_depth = 13, min_samples_leaf = 8, max_features = sqrt, class_weight = balanced, n_jobs = -1, random_state = 42 | Bagging ensemble; min_samples_leaf = 8 prevents leaf-level over-splitting on noisy fault labels. |
| **XGBoost-Full** | n_estimators = 150, max_depth = 6, learning_rate = 0.0407, subsample = 0.925, colsample_bytree = 0.8, tree_method = hist, scale_pos_weight = neg/pos ratio, eval_metric = auc, n_jobs = -1, random_state = 42 | Histogram tree method for speed on 203-dim tabular data; scale_pos_weight aligns class weighting with the other tree models. |
| **LightGBM-Full** | n_estimators = 100, num_leaves = 27, learning_rate = 0.0108, objective = binary, metric = auc, class_weight = balanced, n_jobs = -1, random_state = 42, verbose = -1 | Leaf-wise growth; num_leaves = 27 is well below the 2^max_depth ceiling, intentionally restraining capacity. |
| **CatBoost** | iterations = 200, depth = 7, learning_rate = 0.092, auto_class_weights = Balanced, random_seed = 42, thread_count = -1, verbose = 0, eval_metric = AUC | Symmetric tree structure with native categorical support; auto_class_weights matches the imbalance handling of the other tree baselines. |
| **SVM** | kernel = rbf, C = 5.374, gamma = scale, class_weight = balanced, probability = True, cache_size = 500, random_state = 42 | RBF kernel SVM; probability = True enables Platt-scaled probability output required for AUC. |
| **ARIMA** | (p, d, q) = (1, 0, 0), trend = c, enforce_stationarity = False, enforce_invertibility = False; per-station fitting with rolling forecast window | Univariate AR(1) baseline fitted independently per station; (p, d, q) selected via Optuna over candidate orders {(1,0,0), (1,0,1), (2,0,0), (2,0,1)}. |
| **LSTM-Single** | hidden = 128, num_layers = 1, dropout = 0.389, lr = 2.16e-3, optimizer = AdamW, weight_decay = 1e-4, grad_clip = 2.0, batch_size = 256, max_epoch = 60, patience = 15, pos_weight = neg/pos ratio, random_state = 42 | Single-cell LSTM treated as a unified deep baseline; pos_weight in BCE loss replaces explicit class re-weighting. |
| **GRU-Single** | hidden = 64, num_layers = 1, dropout = 0.362, lr = 1.23e-3, optimizer = AdamW, weight_decay = 1e-4, grad_clip = 2.0, batch_size = 256, max_epoch = 60, patience = 15, pos_weight = neg/pos ratio, random_state = 42 | Standalone GRU at a comparable scale to the GRU expert in FGDSE, but without specialization or gating. |
| **TCN-Standard** | channels = 64, kernel_size = 3 (single branch), dilations = (1, 2, 4, 8), dropout = 0.196, lr = 2.79e-4, optimizer = AdamW, weight_decay = 1e-4, grad_clip = 2.0, batch_size = 256, max_epoch = 80, patience = 15, pos_weight = neg/pos ratio, random_state = 42 | Standard single-branch TCN without multi-scale parallel kernels, SE channel attention, or temporal-attention pooling. |
| **Transformer** | d_model = 128, nhead = 8, num_encoder_layers = 2, dim_feedforward = 256, dropout = 0.397, lr = 2.12e-4, optimizer = AdamW, weight_decay = 1e-4, grad_clip = 2.0, batch_size = 256, max_epoch = 80, patience = 15, pos_weight = neg/pos ratio, random_state = 42 | Encoder-only Transformer with sinusoidal positional encoding; depth chosen from {1, 2, 3} during search. |